\documentclass[11pt]{amsart}
\usepackage{graphicx}
\usepackage{latexsym}
\usepackage{amsfonts}
\usepackage[utf8]{inputenc}
\usepackage{amsthm}
\usepackage{amsmath}
\usepackage{amssymb}
\usepackage[all]{xy} 
\usepackage{mathrsfs} 
\usepackage{stmaryrd} 
\usepackage{eepic}
\usepackage{epic}
\usepackage{eucal}
\usepackage{epsfig}
\usepackage{psfrag}
\usepackage{float}
\usepackage{amssymb}
\usepackage{amsmath}
\usepackage{amsthm}
\usepackage{amsbsy}
\usepackage{bm}
\usepackage{graphicx}
\usepackage{color}
\usepackage[colorlinks=true]{hyperref}
\usepackage{mathrsfs} 
\usepackage{psfrag,epsf}
\usepackage{color}
\usepackage{graphicx}
\usepackage{amsmath}
\usepackage{amsfonts}
\usepackage{enumerate}
\usepackage{latexsym}
\usepackage{amsthm}

\makeatletter\@addtoreset{equation}{section}\makeatother
\makeatletter\@addtoreset{figure}{section}\makeatother
\makeatletter\@addtoreset{table}{section}\makeatother

\newtheorem{theorem}{Theorem}[section]

\newtheorem{example}[theorem]{Example}

\newtheorem{remark}[theorem]{Remark}

\theoremstyle{plain}

\newcommand{\R}{{\mathbb R}}
\newcommand{\C}{{\mathbb C}}
\newcommand{\Z}{{\mathbb Z}}

\newcommand{\T}{{\mathbb T}}

\newcommand{\op}[1]{\!\!\mathop{\rm ~#1}\nolimits}

\newcommand{\scriptop}[1]{\!\!\mathop{\mbox{\rm \scriptsize ~#1}}\nolimits}

\newfont{\gothic}{eufm10 scaled\magstep0}
\newcommand{\got}[1]{\mbox{\gothic #1}}

\title{Symplectic actions of non\--Hamiltonian type}

\begin{document}

\author{\'Alvaro Pelayo}

\begin{abstract}
Hamiltonian symplectic actions of tori on compact 
symplectic manifolds have been extensively studied in the past thirty years, and a
number of classifications have been achieved, for instance in the case that the acting torus is
$n$\--dimensional and the symplectic manifold is $2n$\--dimensional.
In this case the $n$\--dimensional orbits are Lagrangian, so it is  natural to wonder 
whether there are interesting  classes of symplectic actions with Lagrangian orbits, and that are not Hamiltonian. 
It turns out that there are many such classes which contain for example the Kodaira variety, and which can be classified 
in terms of symplectic invariants. The paper reviews several classifications,
 which include symplectic actions having a Lagrangian orbit or a symplectic 
 orbit of maximal dimension.  We make an emphasis on the construction of the 
symplectic invariants, and their computation in examples. 
\end{abstract}

 \dedicatory{In memory of Professor Johannes (Hans) J. Duistermaat (1942--2010)}
\maketitle

\section{Introduction}

The study of Hamiltonian symplectic actions in the past thirty years
 has been a driving force in equivariant symplectic
geometry, and there is an extensive theory, which includes several classification results. 
It is natural to wonder whether there are interesting classes of actions
that are symplectic, but not Hamiltonian, and if so, whether they may be viewed
as part of a larger cohesive framework.  From a geometric view point, being only symplectic is
a natural assumption, and while Hamiltonian actions of maximal dimension appear as symmetries in 
many integrable systems in mechanics, non\--Hamiltonian actions also occur in physics, 
see eg. Novikov~\cite{novikov}, as well as in complex
algebraic geometry, and in topology.

In this paper 
we review the classical results about symplectic {Hamiltonian} actions by Atiyah, 
Guillemin\--Sternberg, and Delzant \cite{atiyah, gs, De}, 
and then present a number of recent results (since around 2002) on the classification of certain classes
of non\--Hamiltonian symplectic torus actions, and which are contained in the work Benoist, Duistermaat, and the author \cite{benoist, benoistcorr, DuPe,Pe}. We outline connections of these works 
with algebraic geometry (in particular Kodaira's classification of complex analytic 
surfaces \cite{kodaira}), with topology (in particular the work of Benson\--Gordon \cite{bg} on torus bundles over tori and nilpotent 
Lie groups), with the work of Guillemin\--Sternberg on multiplicity\--free spaces \cite{multfree}, and 
also with orbifold theory (for instance Thurston's classification of compact $2$\--dimensional
orbifolds). 

Throughout the paper we emphasize the role that the Hamiltonian theory of Atiyah et al. plays in the non\--Hamiltonian theory; essentially, when a non\--Hamiltonian
torus action has a subtorus which acts in a Hamiltonian fashion, the Hamiltonian theory can be applied to symplectically 
characterize the action of such subtori in terms of a combinatorial invariant (a polytope).

The proofs of the results by
Atiyah et al. are covered from a number of view points in the literature (see for
instance Guillemin's book \cite{guillemin}), and
we will not attempt to discuss them further. Instead, we will  describe
in a way accessible to non experts the recent symplectic classifications in~\cite{benoist, benoistcorr, DuPe,Pe}, making an emphasis in the construction of the symplectic invariants that appear in the classification theorems. We also carry out the computation of the invariants in 
examples such as the Kodaira variety \cite{kodaira} (also known as the Kodaira\--Thurston manifold \cite{Th}).

A \emph{symplectic manifold} is a smooth manifold $M$ equipped with a smooth non\--degenerate
closed $2$\--form $\omega$, called a {\em symplectic form}. 
Equivariant symplectic geometry is concerned with the study of
smooth actions of Lie groups $G$ on symplectic manifolds $M$, by means of diffeomorphisms $\varphi \in {\rm Diff}(M)$
which pull\--back the symplectic form $\omega$ to itself $\varphi^*\omega=\omega$ (these are called
\emph{symplectomorphisms}). Actions
satisfying this natural condition are called \emph{symplectic}. In this paper we treat the case when
$G$ is a compact, connected, abelian Lie group, that is, a torus: $T \simeq (S^1)^k$, $k \geq 1$. 
Let $\mathfrak{t}$ be the Lie algebra of $T$, and let $\mathfrak{t}^*$ be its dual Lie algebra.
A fundamental subclass of symplectic
 actions admit what is called al \emph{momentum map}, which is a $\mathfrak{t}^*$\--valued smooth
 function on $M$ which encodes information about $M$ itself, the symplectic form, and the action; 
 such symplectic actions  are called \emph{Hamiltonian}.

The category
of Hamiltonian actions, while large, does not include some simple examples of symplectic actions, for
instance free symplectic actions on compact manifolds, because Hamiltonian
actions on compact manifolds always have fixed points. One striking (and basic) open question
is: {are there non\--Hamiltonian symplectic $S^1$\--actions on compact connected
manifolds with non\--empty \emph{and} discrete fixed point set?} In recent years there has been a flurry of
activity related to this question, see for instance  Godinho~\cite{Go1, Go2}, Jang~\cite{dj0, dj}, Li\--Liu \cite{LiLiu}, 
 Pelayo\--Tolman~ \cite{PeTo}, Sabatini~\cite{Sab}, and Tolman\--Weitsman~\cite{ToWe}; the answer is unknown. 

In the present paper we will not discuss this question because the dimension of the torus is one, and instead 
we focus on symplectic actions of tori of dimension
$k$ on manifolds of dimension $2n$ where $k \geq n$.  When in addition to being symplectic the action is
Hamiltonian, then necessarily  $n=k$, but there are many  non\--Hamiltonian symplectic
actions when $n=k$, and also when $n\geq k+1$.

We will concentrate on smplectic non\--Hamiltonian actions in three distinct cases: when the manifold
is four\--dimensional, when there is an orbit of maximal dimension which is symplectic, 
and when there is a coisotropic orbit (these last two cases are disjoint from each other,
and the four dimensional case overlaps with both of them). The moduli space of symplectic 
actions satisfying one of these conditions
is large, and includes as a particular case Hamiltonian actions of maximal dimension (see \cite{PePiRaSa}
for the description of the moduli space of Hamiltonian actions of maximal dimension on $4$\--manifolds).

Some of the techniques to study Hamiltonian torus actions (see for instance 
the books by Guillemin~\cite{guillemin}, Guillemin\--Sjamaar~\cite{gusj2005}, 
and  Ortega\--Ratiu~\cite{ORbook}) are useful in the study of 
non\--Hamiltonian symplectic torus actions (since many non\--Hamiltonian actions exhibit proper
subgroups which act Hamiltonianly).  In the study of Hamiltonian actions, one tool that  is often used is Morse theory for the (components of the) momentum map of the action. Since there is no momentum map in the classical sense for
a general symplectic action, Morse theory does not appear as a natural tool in the  non\--Hamiltonian case. 
There is an  analogue, however, \emph{circle valued\--Morse theory} (since any symplectic circle action 
admits a circle\--valued momentum map, see McDuff~\cite{MD} and \cite{PeRa}, which is also Morse in a sense) but  it is less 
immediately useful in our setting; for instance a more complicated form of the Morse inequalities holds 
(see Pajitnov~\cite[Chapter 11, Proposition 2.4]{Pajitnov2006} and Farber~\cite[Theorem 2.4]{Farber2004}),
and the theory appears more difficult to apply, at least in the context of non\--Hamiltonian symplectic actions; see \cite[Remark~6]{PeRa} for further discussion in this direction. This could be one reason that non\--Hamiltonian symplectic actions have been studied less in the literature than their Hamiltonian counterparts. However, as this paper
shows, there are rich classes of non\--Hamiltonian symplectic actions, which include  examples of interest not only
in symplectic geometry, but also in algebraic geometry (eg.~the Kodaira variety), differential geometry (eg.~multiplicty free
spaces), and topology (eg.~nilmanifolds over nilpotent Lie groups, orbifold bundles).  
For instance, a compact symplectic manifold endowed with a symplectic action with
Lagrangian orbits (eg. the Kodaira variety) is characterized in terms of four symplectic invariants 
(Theorem~\ref{t2}); if the action is also Hamiltonian, only one invariant remains (the 
polytope $\Delta$).

\vspace{1mm}
\emph{Structure of the paper}. In Section~\ref{ham:sec} we review some of the foundational results on symplectic Hamiltonian torus actions on compact manifolds; the theory is  well documented in the literature, 
so we limit ourselves to the core aspects. In Section~\ref{cos:sec} we  discuss the classification of
 symplectic actions  when there exists a coisotropic orbit.
 In Section~\ref{symp:sec}, we explain the classification when
there exists a symplectic orbit of the same dimension as the acting torus.
In Section~\ref{four:sec} we state the classification of
symplectic actions of $2$\--tori on $4$\--manifolds.
In Section~\ref{proof:sec} we outline the proof of a key result in Section~\ref{cos:sec}. 
Because this paper is introductory we avoid technical statements; references are given throughout for those
 interested in further details.

\section{Symplectic torus actions of Hamiltonian type} \label{ham:sec}

This section treats the case when the action, in addition to  being symplectic, is Hamiltonian. 
Let $(M, \,\omega)$ be  compact, connected, $2n$\--dimensional symplectic manifold.
Let $T$ be a torus.
Suppose that $T$ acts effectively and symplectically on $M$.
Recall the meaning of these notions: a differential $2$\--form
$\omega$ on $M$ is \emph{symplectic} if it is closed, i.e. ${\rm d} \omega=0$, as well as
non\--degenerate. The action $T \times M \to M$ is \emph{effective} if every element in the torus $T$
moves at least one point in $M$, or equivalently
$
\bigcap_{x \in M} \, T_x=\{e\},
$ 
where
$
T_x:=\{t \in T \,\, | \,\, t \cdot x=x\}
$
is the \emph{stabilizer subgroup of the $T$\--action at $x$}.  The action is \emph{free} if
$T_x=\{e\}$ for every $x \in M$. The action $T \times M \to M$
is \emph{symplectic} if $T$ acts by symplectomorphisms, i.e. diffeomorphisms 
$\varphi\colon M \to M$ such that $\varphi^*\omega=\omega.$

A  type of symplectic actions are {\em Hamiltonian} actions.
Let $\got{t}$ be the Lie algebra of $T$ and $\got{t}^*$ its dual.
A symplectic action $T \times M \to M$ is \emph{Hamiltonian} if there
is a smooth map $\mu \colon M \to \mathfrak{t}^*$ such that 
Hamilton's equation
\begin{eqnarray} \label{xts}
-{\rm d} \langle \mu, \, X \rangle =  \textup{i}_{X_M}\omega:=\omega(X_M,\cdot),  \, \, \, \, \, \forall X \in \mathfrak{t},
\end{eqnarray}
holds, where
$
X_M$ is the 
vector field infinitesimal action of $X$ on $M$, 
and the right hand\--side of equation (\ref{xts}) is the one\--form obtained by pairing of the symplectic
form $\omega$ with $X_M$, while the left hand side is the differential of the
real valued function $\langle \mu(\cdot) ,\, X\rangle$ obtained by evaluating elements
of $\got{t}^*$ on $\got{t}$. The Lie algebra $\got{t}_x$ of $T_x$ is the kernel of the linear mapping 
$X\mapsto X_M(x)$ from $\got{t}$ to $\op{T}_x\! M$. 
In the upcoming sections we will use use the notation $\got{t}_M(x):={\rm T}_x(T\cdot x)$, where
$T \cdot x:=\{t \cdot x \, |\, t \in T\}$ is the \emph{$T$\--orbit} that goes through the point $x$.

It follows from equation (\ref{xts}) that Hamiltonian $T$\--actions on compact connected manifolds have fixed points. 
The Atiyah-Guillemin-Sternberg Convexity Theorem (1982, \cite{atiyah,gs}) says that 
$\mu(M)$ is  the convex hull of  the image under $\mu$ of the fixed point set of the $T$\--action.
The polytope $\mu(M)$ is called the \emph{momentum polytope of $M$}. One precedent of this 
result appears in Kostant's article~\cite{kostant}.  Later
Delzant (1988, \cite{De}) proved that if the action if effective and $2\dim T=2n=\dim M$ (in which case the triple
$(M,\omega,T)$ is called a \emph{Delzant manifold} or a \emph{symplectic\--toric manifold}), then 
$\mu(M)$ is a Delzant polytope (i.e. a simple, edge\--rational, and smooth polytope; or in other
words, there are precisely $n$\--codimension $1$ faces meeting at each vertex and
their normal vectors span a $\Z$\--basis of the integral lattice
$T_{\Z}:={\rm ker}\,{\rm exp} \colon {\got t} \to T$). Moreover, Delzant proved that
$\mu(M)$ classifies $(M, \, \omega,\,T)$  in the sense of uniqueness (two Delzant
manifolds are $T$\--equivariantly symplectomorphic if and only if they have the same associated momentum
polytope), and existence (for each abstract Delzant polytope $\Delta$ in $\got{t}^*$ there
exists a Delzant manifold with momentum polytope $\Delta$). The $\dim T$\--dimensional orbits
of a symplectic\--toric manifold are Lagrangian submanifolds, that is, the symplectic form
vanishes along them.

The simplest example of a Hamiltonian torus action is
$
(S^2, \, \omega=\textup{d}\theta \wedge\textup{d}h)
$
equipped with the
rotational circle action $\mathbb{R}/\mathbb{Z}$
about the vertical axis of $S^2$.
This action has
momentum map  $\mu \colon S^2 \to \mathbb{R}$
equal to the height function 
$\mu(\theta, \,h)= h$, and in this case the 
momentum polytope is the interval
$\Delta=[-1,\,1]$.
Another example of a Hamiltonian torus action is
the $n$\--dimensional complex projective space equipped
with a $\lambda$\--multiple, $\lambda>0$, of the Fubini\--Study form
$(\mathbb{CP}^n,\, \lambda \cdot \omega_{\scriptop{FS}})$ and the
rotational $\mathbb{T}^n$\--action induced from the
rotational $\mathbb{T}^n$\--action on the $(2n+1)$\--dimensional
complex plane. This action is Hamiltonian, with
momentum map components given by
$
\mu^{\mathbb{CP}^n,\lambda}_k(z)=\frac{\lambda |z_k|^2}{\sum_{i=0}^n|z_i|^2}$ for each $k=1,\ldots,n$.
The associated momentum 
polytope is
$
\Delta=\textup{convex hull }\{0,\, \lambda {\rm e}_1,\ldots,\lambda {\rm e}_n\},
$
where ${\rm e}_1,\ldots,{\rm e}_n$ are the canonical basis vectors of $\R^n$.
There exists an extensive theory of Hamiltonian actions and related topics,
 see for instance the books by Guillemin~\cite{guillemin}, Guillemin\--Sjamaar~\cite{gusj2005}, 
and Ortega\--Ratiu~\cite{ORbook}.

\section{Symplectic torus actions with coisotropic orbits} \label{cos:sec}

This section is based on parts of \cite{DuPe}; see also Benoist~\cite{benoist, benoistcorr}. 
To make the general theory more transparent, we make an emphasis on examples, the connections with algebraic geometry
and topology, and the construction of the symplectic invariants. Let $(M, \,\omega)$
be a compact, connected, $2n$\--dimensional symplectic manifold.
Let $T$ be a torus.
Suppose that $T$ acts effectively and symplectically on $M$.

\medskip

\paragraph{{\bf Coisotropic orbit condition}} 
\emph{We assume throughout this section that there exists a $T$\--orbit which is a 
coisotropic submanifold.}

\subsection{The  meaning of the coisotropic condition}

The  \emph{principal orbit type of $M$}, denoted by $M_{\scriptop{reg}}$,
is  by definition 
$
M_{\scriptop{reg}}:= \{x \in M \, \, \, \, | \, \, \, \, T_x=\{e\}  \}, 
$
or equivalently, $M_{\scriptop{reg}}$ is the set of points where the $T$\--action
is free. The set $M_{\scriptop{reg}}$ is an open dense subset of $M$, and connected
(since $T$ is connected). An orbit of the $T$\--action is \emph{principal} if it lies
inside of $M_{\scriptop{reg}}$.

A submanifold $C \subset M$ is  \emph{coisotropic} if 
$(\textup{T}_xC)^{\omega_x} \subset \textup{T}_xC$
for all $x \in C$, where $(\textup{T}_xC)^{\omega_x}$
is the symplectic orthogonal complement
to the tangent space to $C$; recall that if
$V$ is a subspace of a symplectic vector space $(W, \, \sigma)$, its
\emph{symplectic orthogonal complement}
$
V^{\sigma}
$
consists of the vectors
$w \in W$ such that $\sigma(w,\,v)=0$ for all $v \in V$. 
A submanifold $C \subset M$ is \emph{Lagrangian} if
$
\omega|_C=0$ and 
$\op{dim}C=\frac{\op{dim}M}{2};$
this is a special case of coisotropic submanifold when the inclusion ``$\subset$" 
is an equality ``$=$".
If $C$ is a coisotropic submanifold of dimension $k$, 
then 
$
2n-k=\op{dim}(\op{T}_x\! C)^{\omega _x}\leq
\op{dim}(\op{T}_x\! C)=k
$
shows that $k\geq n$. The submanifold $C$ has the minimal dimension $n$ 
if and only if $(\op{T}_x\! C)^{\omega _x}=\op{T}_x\! C,$
if and only if $C$ is a Lagrangian submanifold of $M$
(that is, an isotropic submanifold of $M$ of maximal dimension $n$).

It is a consequence of the local normal form (i.e. the symplectic tube theorem~\cite[Section~11]{DuPe}) 
of Benoist \cite{benoist, benoistcorr} and Ortega\--Ratiu
\cite{ortegaratiu} that if a symplectic manifold
admits a torus action with a coisotropic orbit, then this
orbit must be principal. Moreover, the existence of  a coisotropic orbit
in the principal orbit type implies that all orbits in the principal
orbit type are coisotropic. That is,  the existence of
a single coisotropic orbit is equivalent to all principal
orbits being coisotropic. 

\begin{remark}
\normalfont
A symplectic torus action
for which one can show that it has a Lagrangian orbit
falls into the category of actions that we study in this section;
this includes Hamiltonian actions of $n$\--tori (see Section~\ref{ham:sec}). 
\end{remark}

\begin{remark}
\normalfont
Let $f,\, g$ be in the set of $T$\--invariant smooth functions 
$\op{C}^{\infty}(M)^T$ and let $x\in M_{\scriptop{reg}}$.  
Let $\mathcal{H}_f$ be the vector field defined by
$\op{i}_{\mathcal{H}_f}\omega=-{\rm d} f$.
Then
$\mathcal{H}_f(x),\,\mathcal{H}_g(x) \in \got{t}_M(x)^{\omega _x}\cap\got{t}_M(x)$. It follows that 
the Poisson brackets 
$
\{ f,\, g\} :=\omega (\mathcal{H}_f,\,\mathcal{H}_g)
$
 of $f$ and $g$ vanish at $x$.  Since the principal orbit type $M_{\scriptop{reg}}$ is dense 
in $M$, we have that $\{ f,\, g\}\equiv 0$ for all 
$f,\, g\in\op{C}^{\infty}(M)^T$ if the principal orbits 
are coisotropic. Conversely, if we have that $\{ f,\, g\} \equiv 0$ for all 
$f,\, g\in\op{C}^{\infty}(M)^T$, then  $\got{t}_M(x)^{\omega _x}
\subset (\got{t}_M(x)^{\omega _x})^{\omega _x}=\got{t}_M(x)$
 for every  $x\in M_{\scriptop{reg}}$, i.e. $T\cdot x$ is coisotropic. 
Therefore  the principal orbits are coisotropic if and only if $\{ f,\, g\} \equiv 0$  
for all $f,\, g\in\op{C}^{\infty}(M)^T$.\footnote{In Guillemin 
and Sternberg \cite{multfree}, a symplectic manifold with a 
Hamiltonian action of an arbitrary compact Lie group is called 
a {\em multiplicity\--free space} if the Poisson brackets of any pair of 
invariant smooth functions vanish.}
\end{remark}

\begin{remark}
\normalfont
The coisotropic condition is natural in view of the 
local models of symplectic torus actions. 
This observation goes back to Benoist~\cite{benoist}.
\end{remark}

\subsection{Examples} \label{examples}
Many classical examples of symplectic torus actions have coisotropic principal
orbits.

\begin{example} \normalfont \label{ktexample}
(Kodaira variety) The first example of a symplectic torus action with coisotropic principal orbits is
the Kodaira variety \cite{kodaira} (also known as the Kodaira\--Thurston manifold \cite{Th}), which is a torus bundle over a torus
constructed as follows. Consider the product symplectic manifold
$
(\mathbb{R}^2 \times (\mathbb{R}/\mathbb{Z})^2,\,
 \textup{d}x_1 \wedge \textup{d}y_1 +
\textup{d}x_2 \wedge \textup{d}y_2),
$
where $(x_1,\,y_1) \in \R^2$ and $(x_2,\,y_2) \in (\R/\Z)^2$. Consider the action of 
$(j_1, \,j_2) \in \mathbb{Z}^2$ on $(\mathbb{R}/\mathbb{Z})^2$ by the matrix group
consisting of
$$
\left
( \begin{array}{cc}
1 & j_2 \\
0 & 1  \\
\end{array} \right),
$$
where $j_2 \in \mathbb{Z}$. The quotient of this symplectic manifold by  the diagonal action of
$\Z^2$ gives rise to a compact, connected, symplectic $4$\--manifold
\begin{eqnarray} \label{ktm}
({\rm KT},\omega):=(\mathbb{R}^2 \times_{\mathbb{Z}^2} (\mathbb{R}/\mathbb{Z})^2, \,\, \,
 \textup{d}x_1 \wedge \textup{d}y_1 +
\textup{d}x_2 \wedge \textup{d}y_2)
\end{eqnarray}
on which the $2$\--torus
$T:=\mathbb{R}/\mathbb{Z} \times \mathbb{R}/\mathbb{Z}$ acts symplectically and freely, where
the first circle acts on the $x_1$\--component, and the second 
circle acts on the $y_2$\--component (one can check that this action is
well defined). We  denote the Kodaira variety by ${\rm KT}$, and the symplectic
form is assumed. Because the $T$\--action is free, all the orbits are principal, and
because the orbits are obtained by keeping the $x_2$\--component and the
$y_1$\--component fixed, we must have that ${\rm d}x_2={\rm d}y_1=0$,
and hence the aforementioned form $\omega$ vanishes along the
orbits, which are therefore Lagrangian submanifolds of ${\rm KT}$. 
\end{example}

\begin{remark}
\normalfont
The symplectic manifold in formula (\ref{ktm}) fits in  the third case 
in Kodaira \cite[Theorem~19]{kodaira}. Thurston 
rediscovered it \cite{Th}, and observed that there exists no K\"ahler structure on ${\rm KT}$
which is compatible with the symplectic form (by noticing that
the first Betti number $\op{b}_1({\rm KT})$ of ${\rm KT}$ is $3$). It follows that
not only the symplectic action we have described in Example~\ref{ktexample} is not Hamiltonian,
but no other $2$\--torus action on ${\rm KT}$ is Hamiltonian either, since compact
$4$\--manifolds with Hamiltonian $2$\--torus actions 
admit a K\"ahler structure compatible with the given symplectic
form~\cite{guillemin, De}.
\end{remark}

\begin{example}[Non\--free symplectic action]
\normalfont
This is an example of a  non\--Hamiltonian, non\--free symplectic $2$\--torus action on a compact, connected,
symplectic $4$\--manifold. Consider the compact symplectic $4$\--manifold
\begin{eqnarray} \label{r2s2}
(M,\omega):=((\mathbb{R}/\mathbb{Z})^2 \times S^2, \, \, \, \textup{d}x  \wedge \textup{d}y 
+\textup{d}\theta  \wedge \textup{d}h ).
\end{eqnarray}
There is a natural action of the $2$\--torus
$T:=\mathbb{R}/\mathbb{Z} \times \mathbb{R}/\mathbb{Z}$ on expression (\ref{r2s2}),
where the first circle of $T$ acts on the first circle $\R/\Z$ of the left factor of $M$,
and the right circle acts on $S^2$ by rotations (about the
vertical axis); see Figure~\ref{sph}. This $T$\--action is symplectic. 
\begin{figure}[htbp]
  \begin{center}
    \includegraphics[height=2cm, width=7cm]{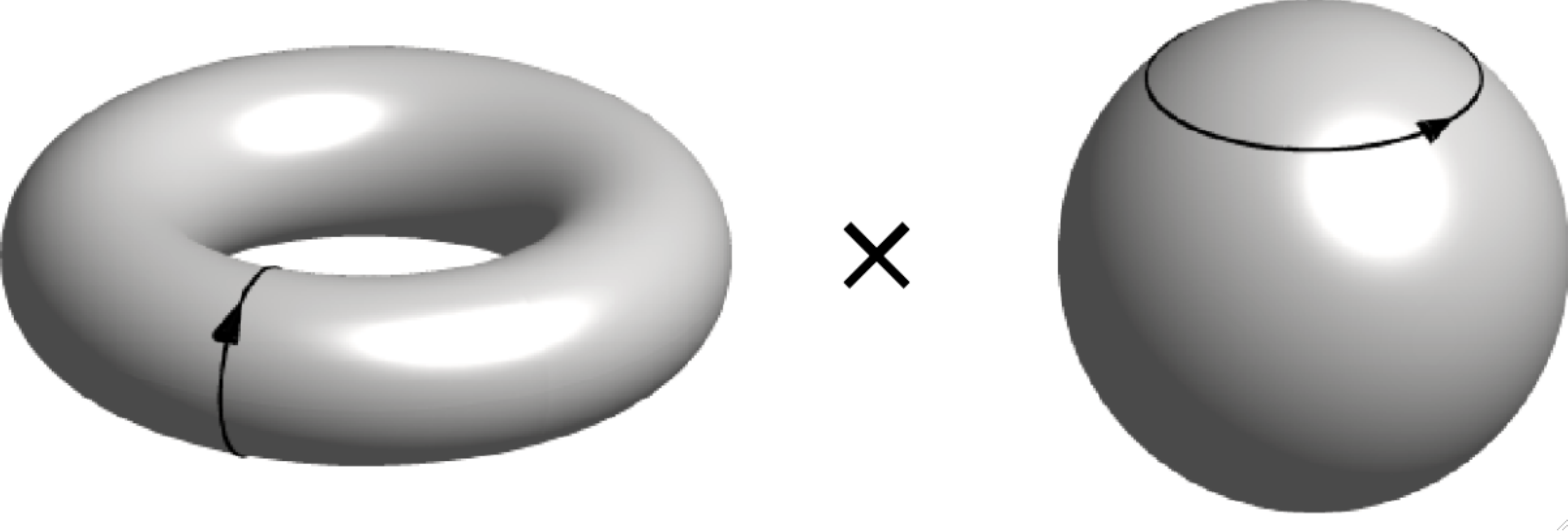}
    \caption{A non\--free symplectic  $2$\--torus action on $S^2 \times (\R/\Z)^2$.}
    \label{sph}
    \end{center}
\end{figure}
However,  it is not a Hamiltonian action because
it does not have fixed points. It is also not
free, because the stabilizer subgroup of a point  $(p,\,q)$, where 
$q$ is the North or South pole of $S^2$, is a circle
subgroup.   In this case the principal orbits are the
products of the circle orbits of the left factor $(\R/\Z)^2$,
and the circle orbits of the right factor (all orbits of the right factor are circles but the
North and South poles, which are fixed points).  Because
these orbits are obtained by keeping the $y$\--coordinate
on the left factor constant, and the height on the right factor
constant, ${\rm d}y={\rm d}h=0$, which implies that
the product form vanishes along the principal orbits, which are
Lagrangian, and hence coisotropic.
\end{example}

\subsection{Enumeration of all examples} \label{t1sec}

The examples in Section~\ref{examples} are described by the following theorem.

\begin{theorem}[\cite{DuPe}] \label{t1}
If a compact, connected, symplectic manifold $M$ admits an effective symplectic 
$T$\--action with a coisotropic orbit\footnote{or Lagrangian orbit}, then
$M$ is the total space of a fibration 
$
F \hookrightarrow M \to B
$ 
with base $B$ being a torus bundle
over a torus, and symplectic toric varietes $F$ as fibers.  On each of the 
toric varietes $F$ a unique subtorus of $T$ acts Hamiltonianly,
any complement of which acts freely by permuting the toric
varietes.
\end{theorem}

\begin{figure}[htbp]
  \begin{center}
    \includegraphics[height=7.5cm, width=7cm]{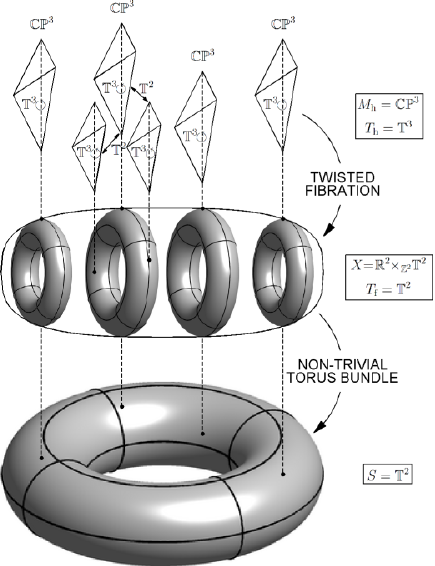}
    \caption{A $10$\--dimensional symplectic manifold with a torus
      action with Lagrangian orbits. The fiber is the toric variety 
      $(\mathbb{CP}^3,\mathbb{T}^3)$, where $\mathbb{T}=\R/\Z$. The base is a non\--trivial
      torus bundle $\R^2 \times_{\Z^2} \mathbb{T}^2$ over the $2$\--torus $\mathbb{T}^2$.}
    \label{LO}
  \end{center}
\end{figure}

Let us describe specifically the model of the symplectic manifold 
$(M,\omega)$ with $T$\--action 
in Theorem \ref{t1}. The symplectic manifold
$(M,\omega)$ is isomorphic (i.e. $T$\--equivariantly symplectomorphic) 
to a fibration  (see Figure~\ref{LO})
$
M_{\scriptop{h}} \hookrightarrow G \times_H M_{\scriptop{h}} \to G/H
$
with the following fiber and base. The fiber $(M_{\scriptop{h}},\,T_{\scriptop{h}})$ is a symplectic\--toric manifold (a toric
variety).  Here $T_{\scriptop{h}}$ is the maximal subtorus of $T$ which acts on
$M$ in a Hamiltonian fashion. The base $G/H$ is a $(T/T_{\scriptop{h}})$\--bundle over a torus $(G/H)/T$,
where  $G$ is (in general possibly) a non\--abelian $2$\--step Nilpotent Lie group defined in terms of the Chern class
of the principal torus bundle $M_{\scriptop{reg}} \to M_{\scriptop{reg}}/T$,
and $H \le G$ is a closed Lie
subgroup of $G$ defined in terms of the holonomy of a certain
connection for the principal torus bundle $M_{\scriptop{reg}}
\to M_{\scriptop{reg}}/T$.

\begin{remark}
\normalfont
Notice that (i) {if the action is free}, then the Hamiltonian subtorus $T_{\scriptop{h}}$
is trivial, and hence $M$ is itself a torus bundle over a torus. Concretely,
$M$ is of the form $G/H$. The Kodaira variety (Example~\ref{ktexample})
is one of these spaces. Since $M$ is a principal torus 
bundle over a torus, it is a nilmanifold for a two\--step 
nilpotent Lie group as explained in  Palais\--Stewart \cite{ps}.  In the
case when this nilpotent Lie group is not abelian,  $M$
does not admit a K\"ahler structure, see Benson\--Gordon \cite{bg}. 
(ii) In the case of  $4$\--dimensional manifolds $M$, item (i) corresponds to the third case 
in Kodaira's description \cite[Theorem~19]{kodaira} of the 
compact complex analytic surfaces which 
have a  holomorphic $(2,\, 0)$\--form that is nowhere 
vanishing, see \cite{DuPecois}\footnote{In~\cite{DuPecois} the authors show
 that  a compact connected symplectic 
$4$\--manifold with a symplectic $2$\--torus action admits an invariant complex structure and give
an identification of those that do not admit a K\"ahler structure with Kodaira's class of complex surfaces
which admit a nowhere vanishing holomorphic $(2, 0)$\--form, but are not a torus or a K3 surface.}. 
As mentioned, these were rediscovered by Thurston \cite{Th}
as the first examples of compact connected symplectic manifolds 
without K\"ahler structure.
(iii) {If on the other hand the action is Hamiltonian}, then
$T_{\scriptop{h}}=T$, and in this case $M$ is itself the toric variety (see~\cite{De, guillemin, DuPeTV} for
the relations between symplectic toric manifolds and toric varieties). 
(iv)
Henceforth,
{we may view the coisotropic orbit case as a twisted mixture of the Hamiltonian case,
and of the free symplectic case}\footnote{The twist affects the manifold topology, 
the torus action, and the sympletic form.}. The twist is
determined by several symplectic invariants of the manifold, which we describe below.
\end{remark}

\subsection{Symplectic ingredients in Theorem \ref{t1}}

Before stating our second theorem on symplectic actions with coisotropic orbits 
we need to understand the construction of
the maximal subtorus $T_{\scriptop{h}} \subset T$, the toric variety 
$M_{\scriptop{h}} \subset M$, as well as how $T$ acts on $M$, 
and how $G$ and $H$ are constructed.

\subsubsection{Construction of $T_{\scriptop{h}}$} \label{thsection}

For this, we first need
to know what our symplectic manifold with $T$\--action looks like locally near 
the $T$\--orbit of a  point $x$. This is a consequence
of the symplectic tube theorem of Benoist \cite{benoist} and Ortega\--Ratiu \cite{ortegaratiu}, which if one
disregards the symplectic structure, says
that there is an isomorphism of Lie groups $\iota$ 
from the stabilizer subgroup $T_x$ onto some $\T ^m$, an open $\T ^m$\--invariant 
neighborhood $E_0$ of the origin in 
$E=(\got{l}/\got{h})^*\times \C ^m$,  
and a $T$\--equivariant diffeomorphism 
$
\Phi \colon K\times E_0 \to U,
$
where $U$ is
an open $T$\--invariant neighborhood of $x$ in $M$ satisfying 
the condition $\Phi (1,\, 0)=x$. Here $K$  is a complementary subtorus of the stabilizer
subgroup $T_x$ in $T$. 
For $t\in T$, let $t_x$ and $t_K$ 
be the unique elements in $T_x$ and $K$, respectively, such that 
$t=t_x\, t_K$. We do not worry about what $\got{l}$ is just yet (it is a particular
vector space), since
the $T$\--action does not affect it as we see next (we give 
a definition of $\got{l}$ in {Section~\ref{Gdef}}). 
The element $t\in T$ acts on $K\times (\got{l}/\got{h})^*\times \C ^m$ 
as
$
t \cdot (k,\,\lambda ,\, z)=(t_K\, k,\,\lambda ,\, \iota (t_x)\cdot z).
$ 
It follows that the stabilizer subgroup of $(k,\,\lambda ,\, z)$ is equal 
to the set of $t_x$ in $T_x$ such that $\iota (t_x)^j=1$ for every 
$j$ such that $z^j\neq 0$. There are 
$2^m$ different stabilizer subgroups $T_y$, $y\in U$. 
Since $M$ is a compact manifold, is follows that there 
are only finitely many different stabilizer subgroups of $T$.  
The product of all the different stabilizer subgroups 
is a subtorus of $T$, which we denote by $T_{\scriptop{h}}$, because  the product of 
finitely many subtori is a compact and connected subgroup of 
$T$, and therefore it is a subtorus of $T$.  One can show that the Hamiltonian subtorus
$T_{\scriptop{h}}$ acts on $M$ in a Hamiltonian fashion, and that
any complementary subtorus $T_{\scriptop{f}}$ to $T_{\scriptop{h}}$
in $T$ must acts freely on $M$.

\subsubsection{Construction of $M_{\scriptop{h}}$} \label{mh}
One can show that
there is a unique smooth distribution $D$ which is integrable and $T$\--invariant,
and the integral manifolds of which are all $(2\,\op{dim}T_{\scriptop{h}})$\--dimensional
symplectic manifolds and isomorphic to each other (i.e. $T_{\rm h}$\--equivariantly
symplectomorphic to each other).
We pick one of them and call it $M_{\scriptop{h}}$.  The symplectic 
form $\omega$ on $M$ restricts to a symplectic form $\omega_{\scriptop{h}}$
on $M_{\scriptop{h}}$. On $(M_{\scriptop{h}},\omega_{\scriptop{h}})$ the
Hamiltonian torus $T_{\scriptop{h}}$ acts Hamiltonianly, and leaving
it invariant.

\subsubsection{Definition of $G$} \label{Gdef}
Using the homotopy formula for the Lie derivative one can show that
there is a unique antisymmetric bilinear form 
$\omega^{\got{t}}$ on $\got{t}$ such that 
\begin{eqnarray} \label{omegax}
\omega _x(X_M(x),\, Y_M(x))=\omega^{\got{t}}(X, Y)
\end{eqnarray}
for every $X,\, Y\in\got{t}$ and every $x\in M$.  
Call $\mathfrak{l}$ the kernel of the antisymmetric bilinear form (the space $\mathfrak{l}$
is for example trivial if the principal orbits are Lagrangian, but this is not
necessarily always the case). Call $N$ the set of linear forms
on $\mathfrak{l}$ which vanish on the Lie algebra $\got{t}_{\scriptop{h}}$ of $T_{\scriptop{h}}$.
Here we are assuming that $\got{t}_{\scriptop{h}} \subset \mathfrak{l}$, which is true because
if $X\in\got{t}_x$, then $X_M(x)=0$.  By slightly abusing the notation, 
we write $N=(\got{l}/\got{t}_{\scriptop{h}})^*$. Set theoretically, we define
$
G = T \times N$. On $G$ we define a non\--standard group structure, in terms of the Chern class of 
$M_{\scriptop{reg}} \to M_{\scriptop{reg}}/T$, which one can naturally identify
with an antisymmetric bilinear form $c \colon N \times N \to \got{l}$. We will
wait to explain this identification to the next section, and we ask the reader
to assume this for the time being, because now we can define the (in general) non\--abelian operation on $G$
\begin{eqnarray} \label{nonabgroup}
(t,\,\zeta )\, (t',\,\zeta ')=
(t\, t'\,\op{e}^{-c(\zeta ,\,\zeta ')/2},\,\zeta +\zeta ').
\end{eqnarray}

\subsubsection{Definition of $H$} \label{hdef:sec}
The subgroup $H$ of $G$ is defined as the set of $(t,\,\zeta )\in G$ such that $\zeta\in P$
and $t\,\tau _{\zeta}\in T_{\scriptop{h}}$. Here $P$ is the so called {\em period lattice} of a natural group action
on  $M/T$ (we explain more in Section~\ref{361}). 
The elements $\tau_{\zeta} \in T$, $\zeta \in P$, encode the holonomy of a certain connection (Step~2 of
Section~\ref{proof:sec}) for the principal torus bundle $M_{\scriptop{reg}} \to M_{\scriptop{reg}}/T$.

\subsubsection{Definition of $G/H$}

The quotient $G/H$ is done with respect to the non standard group structure in expression (\ref{nonabgroup}). 
Moreover, $G/H$ is a nilmanifold.
To go from $G \times_H M_{\scriptop{h}}$ to $G/H$ we cancel the Hamiltonian action
$T_{\scriptop{h}}$, yet we still have an action of the free torus $T/T_{\scriptop{h}} \simeq
T_{\scriptop{f}}$. The quotient $G/H$ is a torus bundle
$
G/H \to (G/H)/(T/T_{\scriptop{h}})
$
over a torus $(G/H)/(T/T_{\scriptop{h}})$, and it follows from work of Palais\--Stewart 
that $G/H$ is a nilmanifold for a $2$\--step nilpotent Lie group. In \cite{DuPe} 
we gave an explicit description of $G/H$ as a nilmanifold.

\subsubsection{The $T$\--action on $M$}
The $T$\--action by translations on the left factor of $G$ passes to an action on 
$G \times_{H} M_{\scriptop{h}}$. We have said that the torus $T$ decomposes
as $T=T_{\scriptop{h}}T_{\scriptop{f}}$, and $T_{\scriptop{h}}$ acts Hamiltonianly 
on each of the fibers, which are identified with the toric variety $(M_{\scriptop{h}},\omega_{\scriptop{h}})$,
leaving them invariants. The free torus $T_{\scriptop{f}}$, which acts freely, simply
permutes the fibers.

\subsubsection{The symplectic form on $M$} 

The following formula bears no weight in the remaining of  this paper. We
write it here to convey that the classification we present is
concrete. Let $\delta a=((\delta t,\,\delta\zeta ),\,\delta x)$, 
and $\delta' a=((\delta 't,\,\delta '\zeta ),\,\delta ' x)$ 
be tangent vectors to the product $G\times M_{\scriptop{h}}$ at 
the point $a=((t,\,\zeta ),\, x)$, 
where we identify each tangent space of  $T$ with 
$\got{t}$. Write 
$X=\delta t+c(\delta\zeta ,\,\zeta )/2$
and 
$
X'=\delta 't+c(\delta '\zeta ,\,\zeta )/2.
$
Let
\begin{eqnarray}
\sigma _a(\delta a,\,\delta 'a)
&=&\omega ^{\got{t}}(\delta t,\,\delta 't)
+\delta\zeta ({X'}_{\got{l}})
-\delta '\zeta (X_{\got{l}})
-\mu (x)(c_{\scriptop{h}}(\delta\zeta ,\,\delta '\zeta ))
\nonumber\\
&&+\, (\omega _{\scriptop{h}})_x(\delta x,\, 
({X'}_{\scriptop{h}})_{M_{\scriptop{h}}}(x))
-(\omega _{\scriptop{h}})_x(\delta 'x,\, 
(X_{\scriptop{h}})_{M_{\scriptop{h}}}(x))
\nonumber\\
&&+\, (\omega _{\scriptop{h}})_x(\delta x,\,\delta 'x).
\label{A*sigma} \nonumber
\end{eqnarray}
Here $X_{\scriptop{h}}$ denotes the 
$\got{t}_{\scriptop{h}}$\--component of $X\in\got{t}$ with respect 
to the decomposition $\got{t}_{\scriptop{h}}\oplus 
\got{t}_{\scriptop{f}}$, where  $\got{t}_{\scriptop{h}}$ is the Lie algebra of
$T_{\scriptop{h}}$, and $\got{t}_{\scriptop{f}}$ is the Lie algebra of $T_{\scriptop{f}}$.  
Similarly holds for $X_{\got{l}}$. The form
$\omega_{\scriptop{h}}$ is the restriction of $\omega$
to the chosen integral manifold $M_{\scriptop{h}}$, cf. Section \ref{mh},
and $c_{\scriptop{h}}(\delta\zeta ,\,\delta '\zeta)$ is the $\got{t}_{\scriptop{h}}$\--component of
$c(\delta\zeta ,\,\delta '\zeta ) \in {\got l}$ in ${\got l}=\got{t}_{\scriptop{h}} \oplus ({\got l} \cap \got{t}_{\scriptop{f}})$.

If $\pi _{M}$ denotes the canonical projection from 
$G\times M_{\scriptop{h}}$ onto 
$G\times _HM_{\scriptop{h}}$, 
then the $T$\--invariant symplectic form 
on $G\times _HM_{\scriptop{h}}$
is the unique two\--form $\beta$ on 
$G\times _HM_{\scriptop{h}}$ such that 
$\sigma ={\pi _{M}}^*\,\beta$.

\subsection{Symplectic classification}
Theorem~\ref{t1} gives a model of the manifold, the torus
action, and the symplectic form, up to isomorphisms. But this is not
a classification, because we do not know how many of these models 
appear, and, for the ones that do appear, we do not know whether
they are isomorphic ($T$\--equivariantly symplectomorphic). 
The first of these is the existence issue,
while the second is the uniqueness issue.

\begin{theorem}[\cite{DuPe}] \label{t2}
The compact, connected, symplectic manifolds $(M, \omega)$ with effective symplectic 
$T$\--actions with some coisotropic orbit are determined (up to $T$\--equivariant
symplectomorphisms) by the following four symplectic invariants: {\rm 1)} The antisymmetric bilinear form
$\omega^{\mathfrak{t}} \colon \mathfrak{t} \times \mathfrak{t} \to
\mathbb{R}$ defined as  the restriction of $\omega$ to the $T$\--orbits; {\rm 2)}
The Hamiltonian torus $T_{\scriptop{h}}$ and its associated momentum
polytope $\Delta$; {\rm 3)} The period lattice $P$ of the maximal the subgroup $N=(\got{l}/\got{t}_{\scriptop{h}})^*$ of $\mathfrak{l}^*$
acting on $M/T$, where $\mathfrak{l}:=\textup{ker}(\omega^{\mathfrak{t}})$; 
{\rm 4a)} The bilinear form $c \colon N \times N \to \mathfrak{l}$ encoding
the Chern class of  the bundle
$M_{\scriptop{reg}} \to M_{\scriptop{reg}}/T$;
{\rm 4b)} An equivalence class $[\tau \colon P \to T]_{\scriptop{exp}(\mathcal{A})} 
\in \op{Hom}_c(P,\, T)/\op{exp}(\mathcal{A})$ encoding the holonomy of a 
certain connection for the bundle $M_{\scriptop{reg}} \to M_{\scriptop{reg}}/T$.
Here $\op{Hom}_c(P,\, T)$ is the set of maps $\tau \colon P \to T$, denoted by
$\zeta \mapsto \tau_{\zeta}$, such that 
$
\tau_{\zeta'}\tau_{\zeta}=\tau_{\zeta+\zeta'}
\op{e}^{c(\zeta',\,\zeta)/2}, 
$
and
$\op{exp}(\mathcal{A})$ eliminates the dependance of the invariants on the choice of base point 
and on the particular choice of connection.
\end{theorem}

For item 1) in Theorem~\ref{t2}  see formula {\rm (\ref{omegax})}. For item 2) see Section{\rm~\ref{thsection}}.
For item 3) see Section{\rm~\ref{hdef:sec}}. For item 4a) see Section{\rm~\ref{Gdef}}.

Section~\ref{thminv} is devoted to a more explicit 
construction of the symplectic invariant 3), and to the construction and 4a), 4b) in Theorem~\ref{t2}. 
While invariants 1) and 2) are relatively straightforward to define, 3) and 4)  are more involved.  

\begin{remark} \label{dez}
\normalfont
(a) Theorem~\ref{t2} is a uniqueness theorem, i.e. we characterize that two manifolds are isomorphic
if and only if they have the same invariants 1) --4). But we do not say which invariants 1) --4) do actually 
appear, which would be the existence theorem. (b) {For simplicity, we have not stated
this existence result here}, but we shall say that, for example, any antisymetric bilinear
form can appear as invariant 1), and any subtorus $S \subset T$ and Delzant polytope
can appear as ingredient 2) etc. This is explained in \cite{DuPe}. 
(c) {Theorem~\ref{t2} is analogous to the
well\--known Delzant theorem} \cite{De}, where Delzant proves that two
manifolds are equivariantly symplectomorphic if and only is they have the same
associated momentum polytope (uniqueness). Then he went on to prove that
one can start from any abstract Delzant polytope and construct
a manifold whose momentum polytope is precisely the
one he started with (existence), and this is the part that we shall
not write here. 
\end{remark}

\subsection{Symplectic invariants in Theorem \ref{t2}} \label{thminv}
We will explain what the invariants in the theorem are.  
Previously we  explained invariants 1) and 2), and we said that later we would explain more
about the natural group action on the orbit space (item 3)), 
the Chern class $c$, and the holonomy invariant (item 4)).

\subsubsection{Explanation of 3): group action on orbit space} \label{361}

Using Leibniz identity
for the Lie derivative one can show thar for each $X\in\got{l}$,
$
\widehat{\omega}(X) := \, -\op{i}_{X_M}\omega
$ 
is a closed 
basic one\--form on $M$.  For each $x\in M$, $\widehat{\omega} (X)_x$ is a linear form 
on $\op{T}_x\! M$ which depends linearly on $X\in\got{l}$, and 
therefore $X\mapsto\widehat{\omega} (X)_x$ is an $\got{l}^*$\--valued 
linear form on $\op{T}_x\! M$, which we denote by 
$\widehat{\omega}_x$. In this way $x\mapsto\widehat{\omega}_x$ 
is an $\got{l}^*$\--valued one\--form on $M$, which we 
denote by $\widehat{\omega}$. With these conventions 
$$
\widehat\omega _x(v)(X)=\widehat{\omega} (X)_x(v)
=\omega _x(v,\, X_M(x)),
\quad x\in M,\; v\in\op{T}_x\! M,\; X\in\got{l}.
$$
Note that the $\got{l}^*$\--valued one\--form 
$\widehat{\omega}$ on $M$ is basic and closed.

In the local model in Section \ref{thsection}
with $x\in M_{\scriptop{reg}}$, 
where $\got{h}=\got{t}_x=\{ 0\}$ and $m=0$, 
at each point the $\got{l}^*$\--valued 
one\--form $\widehat{\omega}$ corresponds to the projection 
$(\delta t,\,\delta\lambda )\mapsto\delta\lambda :
\got{t}\times\got{l}^*\to\got{l}^*$, and 
$\got{t}\times\{ 0\}$ is equal to the tangent space 
of the $T$\--orbit. It follows that   
for every $p\in (M/T)_{\scriptop{reg}}$ the induced 
linear mapping 
$
\widehat{\omega}_p: \op{T}_p(M/T)_{\scriptop{reg}}\to\got{l}^*
$
is a linear isomorphism.  More generally, the strata for the $T$\--action in 
$\overline{M}=K\times (\got{l}/\got{h})^*\times\C ^m$  
are of the form $\overline{M}^J$ in which $J$ is a subset of 
$\{ 1,\,\ldots,\, m\}$ 
and $\overline{M}^J$ is the set of all 
$(k,\,\lambda ,\, z)$ such that $z^j=0$ if and only if $j\in J$. 
If $\Sigma$ is a connected component of 
the orbit type in $M/T$ defined by the subtorus $H$ of $T$ 
with Lie algebra 
$\got{h}$, then for each $p\in\Sigma$ we have 
$\widehat{\omega}_p(X) = 0$ for 
all $X\in {\got h}$, and $\widehat{\omega}_p$ may be viewed as an 
element of $(\got{l}/\got{h})^*={\got h}^0$, 
the set of all linear forms on $\got{l}$ 
which vanish on $\got{h}$.
The linear mapping 
$
\widehat{\omega}_p : 
\op{T}_p\!\Sigma\to (\got{l}/\got{h})^*
$ 
is a linear isomorphism.

Therefore an element $\zeta \in \got{l}^*$ acts naturally on $p \in (M/T)_{\scriptop{reg}}$ by traveling for time
$1$ from $p$ in the direction that $\zeta$ points to. We denote the arrival point
by $p+\zeta$. This action is in general
not well defined at points in $(M/T)\setminus (M/T)_{\scriptop{reg}}$, it is only defined
in the directions of vectors which as linear forms vanish on the stabilizer subgroup
of the preimage under $\pi \colon M \to M/T$. Since $T_{\scriptop{h}}$ is the maximal stabilizer
subgroup, and for each $x$, $T_x \subset T_{\scriptop{h}}$, the additive subgroup
$N:=(\got{l}/\got{t}_{\scriptop{h}})^*$, viewed as the set of linear forms on
$\got{l}$ which vanish on $\got{t}_{\scriptop{h}}$, is the maximal subgroup
of the vector space $\got{l}^*$, which naturally acts on $M/T$. 
The invariant $P$ is the period lattice for the $N$\--action on $M/T$. 
This $\got{l}^*$\--action turns $M/T$ into what is called 
an {\em $\got{l}^*$\--parallel space}. In~\cite[Section ~11]{DuPe} 
one can find a classification of all {\em $V$\--parallel spaces} for any vector space $V$. These
spaces share common characteristics with both locally affine manifolds, 
and manifolds with corners. In~\cite{DuPe} we prove that they are all isomorphic, as $V$\--parallel spaces, to the
product of a closed convex set and a torus. In the case that the $V$\--parallel space
is compact, this convex set is a convex polytope, and in the case of 
the $\got{l}^*$\--parallel space being $M/T$, this
convex polytope is moreover a Delzant polytope. 
A local analysis of the singularities of $M/T$ allows us to define precisely the structure of
$\got{l}^*$\--parallel space. 

Because of the global linear isomorphism
$
\widehat{\omega}_p : 
\op{T}_p\!\Sigma\to (\got{l}/\got{h})^*,
$ 
any 
$\xi\in\got{l}^*$ may be viewed as a constant vector 
field on $(M/T)_{\scriptop{reg}}$.
This is important later in the construction of the connection for the bundle $M_{\scriptop{reg}}
\to M_{\scriptop{reg}}/T$, in terms of which in \cite{DuPe} the authors prove that $M$ is isomorphic to 
$G \times_H M_{\scriptop{h}}$ (the sketch of proof of this is given in Section~\ref{proof:sec}).

\subsubsection{Explanation of 4): Chern class and holonomy invariant} \label{362:sec}
\textup{\,}
\\
\textup{\,}
\emph{The holonomy invariant}.
A vector field $L_{\xi}$ 
on $M_{\scriptop{reg}}$ is a {\em lift of $\xi$}
if  for all $x\in M_{\scriptop{reg}}$ we have that $\op{T}_x\!\pi (L_{\xi}(x))=\xi$. 
Linear assignments of lifts $\xi \in \got{l}^* \mapsto L_{\xi}$ depending linearly on $\xi$
and connections for the principal torus bundle $M_{\scriptop{reg}} \to M_{\scriptop{reg}}/T$
are objects that are equivalent.
Most of the paper~\cite{DuPe} is devoted to the construction
of a  nice connection
\begin{eqnarray} \label{connection}
\xi \in \got{l}^* \to L_{\xi},
\end{eqnarray}
in terms of which one defines the model $G \times_H M_{\scriptop{h}}$, and the
explicit isomorphism (i.e. $T$\--equivariant symplectomorphism) 
between $M$ and this model. By a ``nice" connection
we mean one for which the Lie brackets of two vector fields, and the symplectic
pairings of two vector fields are particularly easy, and zero in most cases. 
Although the vector fields $L_{\xi}$ are smooth on $M_{\scriptop{reg}}$, many of them are singular
(in a severe way, they blow up) on the lower dimensional strata of the orbit type stratification of $M$. 
We shall say more about this 
when we sketch the proof of Theorem~\ref{t2} in Section~\ref{proof:sec}.

In terms of the nice connection (\ref{connection}), 
in \cite{DuPe} we build what we call {\em the holonomy invariant}. We have
to do this because the connection is not unique (although it is essentially unique up to
a factor in the  direction of ${\got t}$), and its holonomy depends on choices
as well.  Indeed, for each $\zeta\in P$ and $p\in M/T$, the curve 
$\gamma_{\zeta}(t):=p+t\,\zeta$,  where
$0\leq t\leq 1$, is a loop in $M/T$. If 
$x\in M$ and $p=\pi (x)$, then the curve 
$
\delta (t)=\op{e}^{t\, L_{\zeta}}(x), 
$
$0\leq t\leq 1$, is called the \emph{horizontal lift of $\gamma _{\zeta}$} 
which starts at $x$, because $\delta (0)=x$, and 
$\delta '(t)=L_{\zeta}(\delta (t))$ is a horizontal 
tangent vector which is mapped by $\op{T}_{\delta (t)}\!\pi$ 
to the constant vector $\zeta$. This implies that  
$\pi (\delta (t))=\gamma _{\zeta}(t)$ for all $0\leq t\leq 1$.  
The element of $T$ which maps 
$\delta (0)=x$ to 
$\delta (1)$ is called the {\em holonomy} 
$\tau _{\zeta}(x)$ of the 
loop $\gamma _{\zeta}$ at $x$
with respect to the 
connection (\ref{connection}). Because $\delta (1)
=\op{e}^{L_{\zeta}}(x)$, we have that
$\tau _{\zeta}(x)\cdot x=\op{e}^{L_{\zeta}}(x)$.
The  element $\tau _{\zeta}(x)$  does depend
on the point $x\in M$, on the period $\zeta\in P$,
or on the choice of nice connection. Hence why we have
to construct a more refined invariant,
the equivalence class of $\tau$ under
a Lie subgroup which we denote by $\op{exp}(\mathcal{A}$),
and which eliminates the dependance on both the choice of connection
and the choice of base point.
\\
\\
\emph{The Chern Class}.  As discussed in Section~\ref{361},  $M/T$ is isomorphic
to the product of a Delzant polytope $\Delta$ and a torus $S$. In fact $S$
is equal to $N/P$, where recall that $N$ is the largest group which
acts naturally on the orbit space $M/T$, and $P$ is its period lattice.
Also, 
$(M/T)_{\scriptop{reg}}\simeq 
\Delta ^{\scriptop{int}}\times (N/P), 
$
where 
$\Delta ^{\scriptop{int}}$ denotes the interior of the Delzant polytope $\Delta$.  
Any connection for the principal $T$\--bundle  $M_{\scriptop{reg}} \to M_{\scriptop{reg}}/T$ has 
a {\em curvature form}, which is a smooth $\got{t}$\--valued 
two\--form on the quotient $M_{\scriptop{reg}}/T$. The cohomology class of 
this curvature form in an element of
$\op{H}^2(M_{\scriptop{reg}}/T,\,\got{t})$, which is independent of the 
choice of the connection.  The $N$\--action on $M/T$ leaves
$M_{\scriptop{reg}}/T\simeq (M/T)_{\scriptop{reg}}$ 
invariant, with orbits isomorphic to the torus $N/P$. The 
pull\--back to the $N$\--orbits defines an isomorphism 
from ${\rm H}^2(M_{\scriptop{reg}}/T,\,\got{t})$ onto 
${\rm H}^2(N/P,\,\got{t})$, which in turn is identified\footnote{as already
observed by \'Elie Cartan.} 
with $(\Lambda ^2N^*)\otimes\got{t}$. 
It follows from the proof in which we construct the nice connection in terms
of which we construct the model of $M$ (cf. Section~\ref{proof:sec} for a sketch), that  
$c \colon N \times N \to \got{l}$, viewed as an element in
$
c\in (\Lambda ^2N^*)\otimes\got{l}\subset(\Lambda ^2N^*)\otimes\got{t},
$ 
equals the negative of the pull\--back to an $N$\--orbit of the 
cohomology class of the curvature form. This shows that 
$c:N\times N\to\got{l}$ is independent of the choice of free
subtorus $T_{\scriptop{f}}$. 
The Chern class $\mathcal{C}$ of the principal $T$\--bundle $\pi :M_{\scriptop{reg}}\to M_{\scriptop{reg}}/T$
is an element of $\op{H}^2(M_{\scriptop{reg}}/T,\, T_{\Z})$.
It is known that the image 
of $\mathcal{C}$ in $\op{H}^2(M_{\scriptop{reg}}/T,\,\got{t})$ 
under the coefficient homomorphism 
$\op{H}^2(M_{\scriptop{reg}}/T,\, T_{\Z})
\to\op{H}^2(M_{\scriptop{reg}}/T,\,\got{t})$ 
is equal to the negative of the cohomology 
class of the curvature form of any connection in the principal 
$T$\--bundle, and hence $c$ represents $\mathcal{C}$.

\subsubsection{Symplectic invariants of the Kodaira variety}

\begin{figure}[htbp] 
  \begin{center}
    \includegraphics[height=5.5cm, width=6cm]{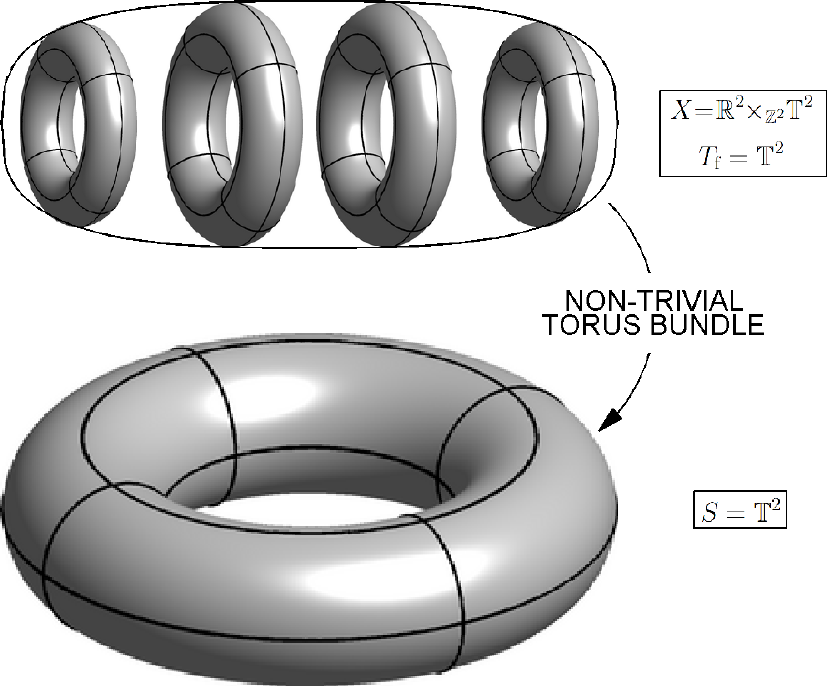}
    \caption{The Kodaira variety endowed with a free symplectic $2$\--torus action with
    Lagrangian orbits.}
       \label{ktPPman}
    \end{center}
    \end{figure}

Recall the Kodaira variety
$M=\R^2 \times_{\Z^2} (\mathbb{R}/\mathbb{Z})^2$ in Example~\ref{ktexample} (see also 
Figure~\ref{ktPPman}). In this case $T=(\R/\Z)^2$, 
$\got{t}\simeq \R^2$ and
its invariants are: 1) the trivial antisymmetric bilinear form: $\omega^{\mathfrak{t}}=0$.
2) The trivial Hamiltonian torus: $T_{\scriptop{h}}=\{[0,\,0]\}$ and the Delzant polytope consisting only of the origin: 
$\Delta=\{(0,\,0)\}$; 3) the period lattice is $P=\mathbb{Z}^2$; 4a)
the bilinear form $c$ representing the Chern class is given by
$
c \colon \mathbb{R}^2 \times \mathbb{R}^2 \to \mathbb{R}^2$,
defined by $c(e_1,\,e_2)=e_1.
$
4b) The holonomy invariant is the class of 
$\tau$
given by
$
\tau_{e_1}=\tau_{e_2}=[0,\,0]$.
In this case
$G=(\mathbb{R}/\mathbb{Z})^2 \times \mathbb{R}^2$,
$M_{\scriptop{h}}=\{p\}$, and
$H=\{[0,\,0]\} \times \mathbb{Z}^2$. The model of $M$ is 
$
G \times_H M_{\scriptop{h}} \,\, \simeq \, \,
G/H \,\, \simeq \, \mathbb{R}^2 \times_{\mathbb{Z}^2} (\mathbb{R}/\mathbb{Z})^2.
$

\section{Symplectic torus actions with maximal symplectic orbits} \label{symp:sec}

This section is based on parts of \cite{Pe}. Let $(M, \,\omega)$
be a compact, connected, $2n$\--dimensional symplectic manifold.
Let $T$ be a torus.
Suppose that $T$ acts effectively and symplectically on $M$.

\medskip

\paragraph{{\bf Symplectic orbit condition}} \emph{Throughout this section we 
assume that there exists a $\op{dim}T$\--dimensional orbit which is a symplectic submanifold.}

\subsection{The meaning of the symplectic condition}

A submanifold 
$C$ of the symplectic manifold  $(M,\omega)$ is \emph{symplectic} if the restriction $\omega|_C$ 
of the symplectic form $\omega$ to $C$ is symplectic. Let $\mathfrak{t}_x$ denote the Lie algebra of 
$T_x$. In Section \ref{Gdef}
we observed that $\mathfrak{t}_x \subset \textup{ker}\, \omega^{\mathfrak{t}}$. Since 
the antisymmetric bilinear form $\omega^{\got{t}}$ introduced therein
(which gives the restriction of $\omega$ to the $T$\--orbits)
is non\--degenerate because of the symplectic condition, 
its kernel $\textup{ker}\, \omega^{\mathfrak{t}}$ is trivial. Hence $\mathfrak{t}_x$ is the trivial vector space, and hence $T_x$, which is
a closed and hence compact subgroup of $T$, must be a finite group (which is of course
abelian since $T$ is abelian). Since $T$ is compact, the action of $T$ on 
$M$ is a proper action, the mapping
$
t \mapsto t \cdot x \colon T/T_x \to T \cdot x
$
is a diffeomorphism, and, in particular, the dimension of the quotient group $T/T_x$ equals the dimension 
of $T \cdot x$. Since each $T_x$ is finite, the dimension of $T/T_x$ equals $\op{dim}T$, and
hence every $T$\--orbit is $\op{dim}T$\--dimensional. Since
the symplectic form $\omega$ restricted to any $T$\--orbit of the $T$\--action 
is non\--degenerate,
$T\cdot x$ is a symplectic submanifold of $(M, \, \omega)$.
Therefore, we can conclude that
at least one $T$\--orbit is a $\op{dim}T$\--dimensional 
symplectic submanifold of the symplectic manifold $(M, \, \omega)$ if and only if
every $T$\--orbit is a $\op{dim}T$\--dimensional 
symplectic submanifold of $(M, \, \omega)$.

Then there exists only finitely many different subgroups of $T$  which occur
as stabilizer subgroups of the action of $T$ on $M$, and each of them is a finite group.
Indeed, we know that every stabilizer subgroup of the
action of $T$ on $M$ is a finite group. 
It follows from  the tube theorem of Koszul (cf. \cite{koszul} or \cite[Theorem~2.4.1]{DuKo}),
that in the case of a compact smooth manifold equipped
with an effective action of a torus $T$, there exists only 
finitely many different subgroups of $T$ which occur
as stabilizer subgroups.

\begin{remark}\normalfont
In the case when $\dim M=4$ and $\dim T=2$,   
the stabilizer subgroup of the $T$\--action at every point in $M$ is
a cyclic abelian group (cf. \cite[Lemma 2.2.6]{Pe}). This follows from
the fact that a finite group acting linearly on a disk must
be a cyclic group acting by rotations (after application of the symplectic tube theorem).
\end{remark}

\subsection{Examples}

\begin{figure}[htbp]
  \begin{center}
    \includegraphics[height=3.4cm, width=6.4cm]{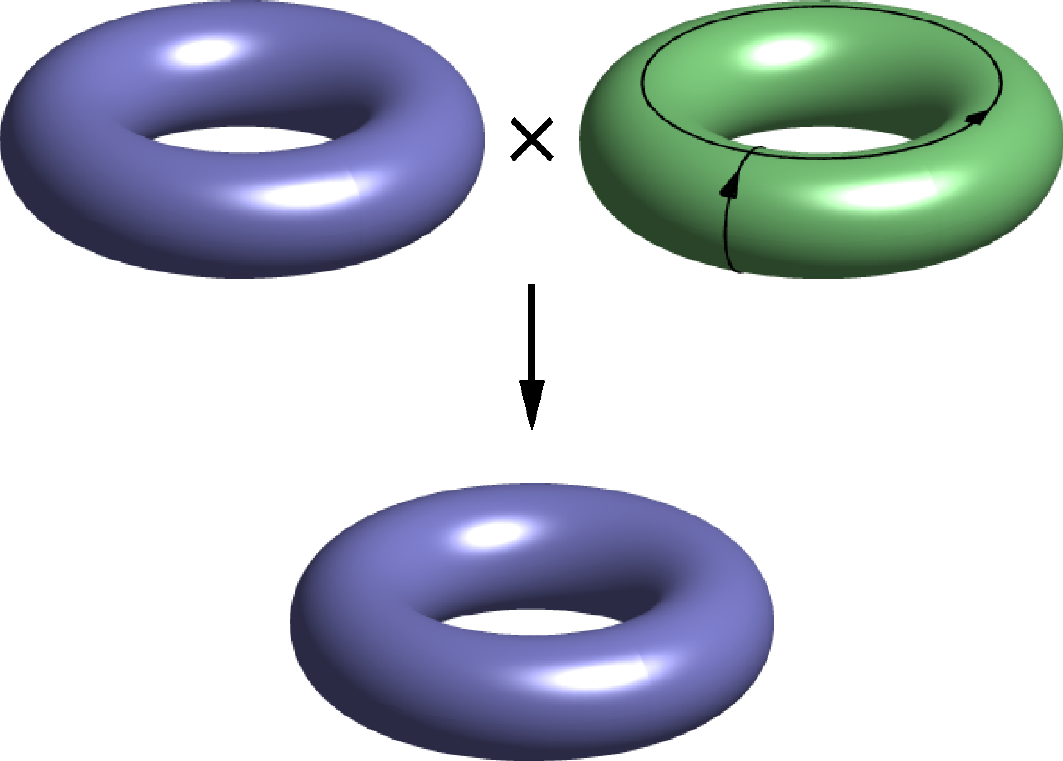}
    \caption{A $4$\--dimensional  torus $(\R/\Z)^2 \times (\R/\Z)^2$ endowed with a symplectic $2$\--torus action with symplectic
    orbits on the right factor. The orbit space is $(\R/\Z)^2$.}
     \label{torustorus}
    \end{center}
\end{figure}

\begin{example}[Free symplectic action]
\normalfont
The $4$\--dimensional torus
$
M:=(\R/\Z)^2 \times (\R/\Z)^2
$
endowed with the
standard symplectic form, on which the $2$\--dimensional torus $(\R/\Z)^2$
acts by multiplications on two of the copies of $\R/\Z$ inside
of $(\R/\Z)^4$, is symplectic manifold with symplectic orbits (see Figure~\ref{torustorus}). 
The $T$\--orbits are symplectic $2$\--tori, and $M/T=(\R/\Z)^2$.
\end{example}

\begin{example}[Free symplectic action]
\normalfont
\label{ex1}
Let $(M, \, \omega):=S^2 \times (\R/\Z)^2$
endowed with the product symplectic form (of the standard symplectic (area) form
on $(\R/\Z)^2$ and the standard area form on $S^2$). 
\begin{figure}[htbp]
  \begin{center}
    \includegraphics[height=1.9cm, width=6.35cm]{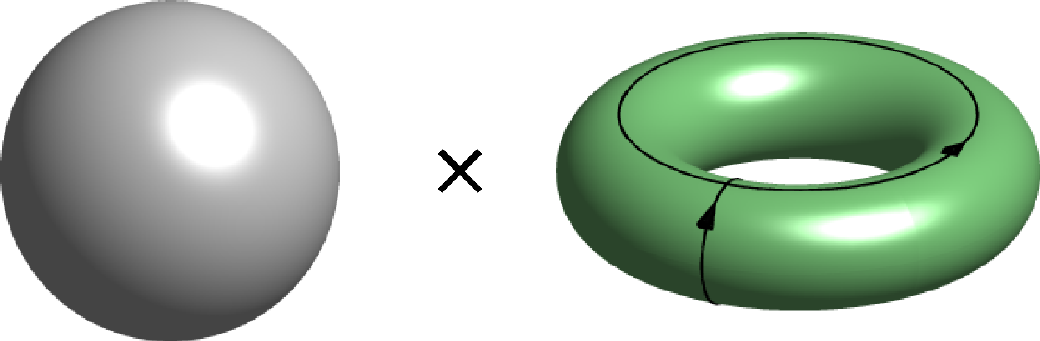}
    \caption{A free symplectic action of a $2$\--torus on the right factor of the product $S^2 \times (\R/\Z)^2$.}
    \label{ttorus00}
    \end{center}
\end{figure}
Let $T:=(\R/\Z)^2$ act on $M$ by translations on the right factor.  This is a free
action on $M$ the orbits of which are symplectic $2$\--tori. In this case
$M/T=S^2$. 
\end{example}

\begin{example}[Non\--free symplectic action] \normalfont
 \label{ex2}
Let $P:=S^2 \times (\R/\Z)^2$ equipped with the product symplectic form of the standard symplectic (area) form
on  $S^2$ and the standard area form on the sphere $(\R/\Z)^2$. The $2$\--torus $T:=(\R/\Z)^2$ acts freely by translations on the right factor
of $P$, see Figure~\ref{ttorus00}.
Let the finite group $\Z/2\,\Z$ act on $S^2$ by rotating
each point horizontally by $180$ degrees, and let  $\Z/2\,\Z$ act
on $(\R/\Z)^2$ by the antipodal action on the first circle $\R/\Z$. 
The diagonal action of $\Z/2\,\Z$ on $P$ is free. Therefore, the quotient space 
$
S^2 \times_{\Z/2\,\Z} (\R/\Z)^2
$ 
is a smooth manifold.
Let $M:=S^2 \times_{\Z/2\,\Z} (\R/\Z)^2$ be endowed 
with the symplectic form $\omega$ and $T$\--action inherited from
the ones given in the product $S^2 \times (\R/\Z)^2$, where $T=(\R/\Z)^2$. The action
of $T$ on $M$ is not free, and the $T$\--orbits are symplectic $2$\--dimensional tori.
The orbit space $M/T$ is $S^2/(\Z/ 2\, \Z)$, which is a 
smooth orbifold with two singular points of order $2$, the South and North poles of $S^2$. 
\end{example}

\subsection{Enumeration of all examples}

The following provides an enumeration of all examples that can appear.

\begin{theorem}[\cite{Pe}] \label{t3}
If a compact, connected, symplectic manifold $M$ admits an effective symplectic 
$T$\--action with a $\op{dim}T$\--dimensional symplectic orbit, then 
$M$ is isomorphic (i.e. $T$\--equivariantly symplectomorphic) to an 
orbifold bundle $\widetilde{M/T} \times_{\Gamma} T$
where  $\widetilde{M/T}$ is the orbifold universal cover\footnote{The universal cover and the 
fundamental group must be based at the same point $p_0$.} of $M/T$ and $\Gamma:=\pi^{\textup{orb}}_1(M/T)$, 
$\Gamma$ acts on $T$ 
by means of the monodromy $\mu \colon \Gamma \to T$ of the flat connection
 $\{(\textup{T}_x(T \cdot x))^{\omega_x}\}_{x \in M}$ of symplectic orthogonal complements to
 the tangent spaces to the $T$\--orbits, and 
 $\Gamma$ acts on the Cartesian product $\widetilde{M/T} 
\times T$ by the diagonal action $x \, (y,\,t)=(x \star y^{-1},\,\mu(x) \cdot t)$, where 
$\star \colon \Gamma \times \widetilde{M/T} \to \widetilde{M/T}$ denotes the natural
action of $\Gamma$ on $\widetilde{M/T}$. The $T$\--action and symplectic form on $\widetilde{M/T} \times_{\Gamma}T$ are inherited from
the $T$\--action on the right factor of $\widetilde{M/T} \times T$ and the natural product form.
\end{theorem}

\begin{remark}
\normalfont
The article \cite{DuPesymp} shows that the first Betti number of $M/T$ is equal to the first
Betti number of $M$ minus the dimension of $T$.
\end{remark}

\subsection{Ingredients of Theorem \ref{t3}}

\subsubsection{$M/T$ is a smooth orbifold}

We denote by 
$
\pi \colon M \to M/T
$ 
the canonical projection $\pi(x):=T \cdot x$. The orbit space $M/T$ is endowed with the maximal topology for which $\pi$ is continuous (which is a Hausdorff topology). Since $M$ is compact and connected, $M/T$ is compact and
 connected. If $C$ is a connected component of an orbit type 
 $
 M^H:=\{x \in M \, |\, T_x=H\},
 $ 
 the action of $T$ on $C$ induces a 
free and proper action of the quotient torus $T/H$ on $C$, and the image 
$\pi(C)$ 
has a unique smooth manifold structure for which
$\pi  \colon C\to\pi(C)$ is a principal $T/H$\--bundle. 

\begin{remark}
\normalfont
In general, $M/T$ is not a smooth manifold, cf. Example \ref{ex2}. 
\end{remark}

Our next goal is to explain how $M/T$ has a natural smooth orbifold structure. 
By the tube theorem (see for instance~\cite[Theorem~2.4.1]{DuKo}) if $x \in M$ there is a $T$\--invariant open neighborhood $U_x$ of 
 $T \cdot x$ and a $T$\--equivariant diffeomorphism  $\Phi_x \colon U_x \to T \times_{T_x} D_x$, 
 where $D_x$ is an open disk centered at the origin in
 $\C^{k/2}$, $k:=\op{dim}M-\op{dim}T$. Here the stabilizer $T_x$ acts by linear transformations on the disk $D_x$.
The action of $T$ on $T \times_{T_x} D_x$ is induced by the translational action of $T$ on
the left factor of the product $T \times D_x$. 
The $T$\--equivariant diffeomorphism $\Phi_x$ induces a homeomorphism  
$\widehat{\Phi}_x \colon D_x/T_x \to \pi(U_x)$, and we have a commutative diagram 
\begin{eqnarray} \label{keydiagram}
\xymatrix{ 
T \times D_x \ar[r]^{\pi_x}      &  T  \times_{T_x} D_x \ar[d]^{p_x}  \ar[r]^{\Phi_x}  & U_x  \ar[d]^{\pi|_{U_x}} \\
 D_x \ar[u]^{i_x} \ar[r]^{\pi'_x} &  D_x/T_x   \ar[r]^{\widehat{\Phi}_x} &      \pi(U_x)} 
\end{eqnarray} 
In this diagram $\pi_x$, $\pi'_x$, $p_x$ are the canonical projections and $i_x$ is the inclusion.
Let $\phi_x:= \widehat{\Phi}_x \circ \pi'_x$. In~\cite{Pe} it is shown that the collection
 $
\widehat{\mathcal{A}}:=\{(\pi(U_x), \, D_x, \, \phi_x, \, T_x)\}_{x \in M}$ 
is an orbifold atlas for  $M/T$.  We call $\mathcal{A}$ the class of atlases
 equivalent to  $\widehat{\mathcal{A}}$.
 We simply write $M/T$ for $M/T$ endowed with 
the class $\mathcal{A}$. 

\subsubsection{The connection $\Omega:=\{(\textup{T}_x(T \cdot x))^{\omega_x}\}_{x \in M}$ is flat}

The collection of subspaces $\Omega:=\{\Omega_x\}_{x \in M}$
where $\Omega_x:=(\op{T}_x(T \cdot x))^{\omega_x}$, is a smooth distribution on $M$ and 
$\pi \colon M \to M/T$ is a smooth principal $T$\--bundle
of which $\Omega$ is a $T$\--invariant flat connection.
To prove this one uses the diagram (\ref{keydiagram}), 
the construction of the orbifold structure on the orbit space $M/T$, and
the symplectic tube theorem  (cf. Benoist \cite[Prop.\,1.9]{benoist} and Ortega and Ratiu \cite{ortegaratiu}).
Put it differently, $\Omega$ is a 
smooth $T$\--invariant $(\op{dim}M-\op{dim}T)$\--dimensional foliation of $M$.

\subsubsection{$M/T$ is a good orbifold, and hence $\widetilde{M/T}$ is a manifold}

Let
$\mathcal{I}_x$ be the maximal integral manifold of $\Omega$. The inclusion  $i_x \colon \mathcal{I}_x \to M$
is an injective immersion and $\pi \circ i_x \colon \mathcal{I}_x \to M/T$
is an orbifold covering . Since $\widetilde{M/T}$ covers any covering of $M/T$, 
 it covers $\mathcal{I}_x$, which is a manifold. Because a covering
of a manifold is itself a manifold, the universal cover $\widetilde{M/T}$ is a manifold.
Readers unfamiliar with orbifolds may  consult \cite[Section~9]{Pe}.

\subsubsection{The orbifold fundamental group $\Gamma$} \label{444:sec}
Let $\mathcal{O}$ be a connected (smooth) orbifold; the case we are interested in is
$\mathcal{O}=M/T$, but the following notions hold in general. The \emph{orbifold
fundamental group $\pi^{\textup{orb}}_1(\mathcal{O},
\, x_0)$ of $\mathcal{O}$ based at $x_0$}
is the set of homotopy classes of orbifold loops that have the same initial
and end point $x_0$. The operation on this set is the classical composition law 
by concatenation of loops. The set $\pi^{\textup{orb}}_1(\mathcal{O},
\, x_0)$ endowed with this operation is a group. It is an exercise to check (as in the classical case) that 
changing the base point results in an isomorphic group.  Thurston proved that any connected 
(smooth) orbifold $\mathcal{O}$ posseses an orbifold covering\footnote{a unique up to equivalence.} 
$\widetilde{\mathcal{O}}$ which is universal (i.e. it covers any other cover), and the 
orbifold fundamental group of which  based at any regular point is trivial.

\subsubsection{$T$\--action on $M$}

The $T$\--action on $\widetilde{M/T} \times_{\Gamma}T$ 
 is the $T$\-- action inherited 
 from the action of $T$ by translations on the right factor of $\widetilde{M/T} \times T$.

\subsubsection{Symplectic form on $M$}

The orbit space $M/T$, which we have seen is a smooth orbifold, comes endowed with a symplectic structure. 
In order to define it, we review the notion of symplectic form on a smooth orbifold $\mathcal{O}$. 
A \emph{(smooth) differential form $\sigma$} (respectively a \emph{symplectic form}) on $\mathcal{O}$
is a collection $\{\widetilde{\sigma}_i\}$, in which $\widetilde{\sigma}_i$ is a 
$\Gamma_i$\--invariant differential form (respectively a symplectic form) on each
$\widetilde{U}_i$, and such that any two of these forms coincide on overlaps (by what we mean that
for each $x \in \widetilde{U}_i$, $y \in \widetilde{U}_j$ with $\phi_i(x)=\phi_j(y)$, there is a diffeomorphism
$\psi$ from a neighborhood of $x$ to a neighborhood of $y$
such that  $\phi_j\circ \psi=\phi_i$ on it, and $\psi^*\sigma_j=\sigma_i$). 
A \emph{symplectic orbifold} is a pair consisting of smooth orbifold and a symplectic form on it.

Let $\sigma$ be a differential form on an orbifold $\mathcal{O}$ given 
by $\{\widetilde{\sigma}_i\}$, where each $\widetilde{\sigma}_i$ is a 
$\Gamma_i$\--invariant differential form on $\widetilde{U}_i$. 
If $f \colon \mathcal{O}' \to \mathcal{O}$ is an orbifold diffeomorphism, the 
\emph{pull\--back} of $\sigma$ is the 
unique differential form $\sigma'$ on $\mathcal{O}'$ given by $\widetilde{f}^*\widetilde{\sigma}_i$
on each  $\widetilde{f}^{-1}(\widetilde{U}_i)$. We shall use the notation $\sigma':=f^*\sigma$.
Similarly one defines the pullback under a principal $T$\--orbibundle $p \colon Y \to X/T$ of $\sigma$. 
The symplectic orbifolds
$(\mathcal{O}_1,\nu_1), \, (\mathcal{O}_2, \nu_2)$ are \emph{symplectomorphic}
if there exists an orbifold diffeomorphism $f \colon \mathcal{O}_1 \to \mathcal{O}_2$
with $f^*\nu_2=\nu_1$ (analogously to  the case of manifolds, the map $f$ is called
an \emph{orbifold symplectomorphism}).

In our case of $\mathcal{O}=M/T$, one can show that there exists a unique $2$\--form $\nu$ on  $M/T$ such that
$\pi^* \nu|_{\Omega_x}=\omega|_{\Omega_x}$ for every $x \in M$.
Moreover, $\nu$ is a symplectic form. Therefore $(M/T,\, \nu)$
is a compact, connected, symplectic orbifold.  The symplectic form on $\widetilde{M/T}$ is the pullback by the 
covering map $\widetilde{M/T} \to M/T$ of  $\nu$ and the symplectic form on 
$T$ is the unique $T$\--invariant symplectic form determined by $\omega^{\mathfrak{t}}$.
The symplectic form on $\widetilde{M/T} \times T$ is the product symplectic form. 
The symplectic form on $\widetilde{M/T}\times_{\Gamma}T$ 
is induced on the quotient by the product form.

\subsection{Symplectic classification}

The following is the classification theorem, up to $T$\--equivariant symplectomorphisms. While in
Theorem~\ref{t3} we make no assumption on the dimension of $T$, in order to give a classification
we need to assume that $\dim T=\dim M-2$. The reason is that in this case the orbit space
$M/T$ is $2$\--dimensional and $2$\--dimensional orbifolds are classified (by Thurston), unlike their
higher dimensional analogues.

\begin{theorem}[\cite{Pe}] \label{t4}
The compact, connected, symplectic $2n$\--dimensional manifolds $(M,\omega)$ equipped
with an effective symplectic $T$\--action with some $(2n-2)$\--dimensional symplectic $T$\--orbit
are classified by the following four symplectic invariants: {\rm 1)} The non\--degenerate  antisymmetric 
bilinear  $\omega^{\mathfrak{t}} \colon \mathfrak{t} \times \mathfrak{t} \to
\mathbb{R}$: the restriction of the symplectic form $\omega$ to the $T$\--orbits; {\rm 2)} 
The Fuchsian signature $(g;\, \vec{o})$ encoding the smooth orbisurface $M/T$; {\rm 3)} 
The symplectic area $\lambda$ of $M/T$; {\rm 4)} The monodromy invariant of the connection $\Omega$ of symplectic orthocomplements to the tangent
spaces to the $T$\--orbits:
$$
 \mathcal{G}_{(g, \, \vec{o})} \cdot 
 ((\mu_{\textup{h}}(\alpha_i), \, \mu_{\textup{h}}(\beta_i))_{i=1}^{g}, \
 (\mu_{\textup{h}}(\gamma_k))_{k=1}^m) \in T^{2g+m}/\mathcal{G}_{(g, \, \vec{o})},
$$
where
$$
\mathcal{G}_{(g; \, \vec{o})}:=\Big\{   
\left( \begin{array}{cc}
A & \textup{0}  \\
C & D 
\end{array} \right) 
 \, \, \,  | \, \, \, A \in \textup{Sp}(2g, \, \Z), \, \, \, D \in \mathcal{M}\textup{S}^{\vec{o}}_m \Big\},
 \nonumber
$$
the map $\mu_{\textup{h}} \colon \textup{H}_1^{\textup{orb}}(M/T, \, \mathbb{Z}) \to T$ is the monodromy homomorphism of $\Omega$, 
the collection $\{\alpha_i,\beta_i\}  \subset \textup{H}_1^{\textup{orb}}(M/T, \, \mathbb{Z})$ 
is a symplectic basis of a maximal free subgroup, and 
$\{\gamma_k\} \subset \textup{H}_1^{\textup{orb}}(M/T, \, \mathbb{Z})$ is a geometric torsion basis.
\end{theorem}

\normalfont

\begin{remark}
\normalfont
Theorem~\ref{t4} is a uniqueness theorem, i.e. we characterize that two manifolds are isomorphic
($T$\--equivariantly symplectomorphic) if and only if they have the same symplectic invariants 1)\--4). 
Theorem~\ref{t4} 
does not say which symplectic invariants 1)\--4) appear; such a statement would be given
by an existence theorem.  In order to keep the presentation simple and the paper succinct, we do not state
the aforementioned existence result here, but we shall say that, for instance, any non\--degenerate antisymetric bilinear form can appear as symplectic invariant 1), and any tuple $(g;\, \vec{o})$ 
such that  
$(g; \, \vec{o})$ is not of the form $(0;\,o_1)$ or of the form $(0; \, o_1, \, o_2)$ with
$o_1 < o_2$ also shows up. This is explained in \cite{Pe}.
\end{remark}

\begin{remark}
\normalfont
Theorem~\ref{t4}  is analogous to the uniqueness part of Delzant's Theorem (\cite{De}) which states that two
Delzant manifolds (also known as symplectic\--toric manifolds) are isomorphic (i.e. $T$\--equivariantly
symplectomorphic) if and only is they have the same
associated momentum polytope (uniqueness), which is a particular type of polytope known as
a \emph{Delzant polytope}\footnote{Simple, rational, and smooth polytope}. Then Delzant went on to proving that
one can start from any abstract Delzant polytope and construct
a symplectic manifold with a torus action whose momentum polytope is precisely the
one he started with (existence); and this is the part that we shall
not write here. 
\end{remark}

\subsection{Symplectic invariants in Theorem \ref{t4}}

Ingredient 1) was explained in Section~\ref{cos:sec}. Ingredients 2) and 3) are
easy to explain, which we do briefly below. We focus on describing  ingredient 4).

\subsubsection{Invariant 2): the orbifold invariant}

Let $\mathcal{O}$ be a smooth orbisurface with $m$ cone points $p_1, \ldots, p_m$.
The \emph{Fuchsian signature of $\mathcal{O}$} is $(g; \, \vec{o})$
where
 $g$ stands for the  genus of the surface underlying $\mathcal{O}$,
$o_k>0$ is the order of $p_k$, and $\vec{o}=(o_1,\ldots,o_m)$, where
$o_k \le o_{k+1}$ for all $1 \le k \le m-1$. 
We denote by $\textup{sig}(\mathcal{O})$ the \emph{Fuchsian signature of $\mathcal{O}$}. 
The orbit space $M/T$, which is a smooth orbisurface, is topologically classified by $(g; \, \vec{o})$.

\subsubsection{Invariant 3): the area invariant}

The symbol $\int_{\mathcal{O}} \sigma$ denotes 
the integral of a differential form $\sigma$ on the orbifold
$\mathcal{O}$. If $(\mathcal{O},\, \sigma)$
is a symplectic orbifold, $\int_{\mathcal{O}} \sigma$
is the \emph{total symplectic area of $(\mathcal{O}, \,
\sigma)$}. By the orbifold Moser's theorem \cite[Theorem~3.3]{MW}, 
if $(M, \, \omega)$, $(M', \, \omega')$ are  compact, connected, symplectic manifolds equipped
with an action of a torus $T$ of dimension $\dim T=\dim M-2$, for which the 
$T$\--orbits are $\op{dim}T$\--dimensional symplectic
submanifolds of $(M, \, \omega)$, and such that
$\int_{M/T} \nu=\int_{M'/T} \nu'$ and  ${\rm sig}(M/T)={\rm sig}(M'/T)$,
then $(M/T, \, \nu),(M'/T,\, \nu')$ are symplectomorphic.

\subsubsection{Invariant 4): the monodromy invariant}

Let $g, \, n, \, o_k $, $1 \le k \le m$, be non\--negative integers. Let $\Sigma$ be a 
compact, connected, orientable smooth orbisurface,  with
underlying topological space a surface of genus $g$ with $m$ singular points
$p_1,\ldots,p_m$ of orders $o_k$. The orbifold fundamental group  is
$
\pi_1^{\textup{orb}}(\Sigma, \, p_0)\!=\!\langle 
\{\alpha_i,\, \beta_i\}_{i=1}^g,   \{\gamma_k\}_{k=1}^m  | 
\prod_{k=1}^m \gamma_k=\prod_{i=1}^g [\alpha_i, \, \beta_i], \,  \gamma_k^{o_k}=1, \, 1 \le k \le m
\rangle,$
where $\alpha_i, \, \beta_i$, $1 \le i \le g$ is
a symplectic basis of the surface underlying $\Sigma$, and the $\gamma_k$ are the homotopy classes of small loops
$\widetilde{\gamma}_k$ around the orbifold singular points. 
By abelianizing $\pi_1^{\textup{orb}}(\Sigma, \, p_0)$ we get
 $
\op{H}_1^{\textup{orb}}(\Sigma,
\, \Z)\!=\!\langle 
\{\alpha_i,\, \beta_i\}_{i=1}^g,  \, \{\gamma_k\}_{k=1}^m \, \, | \, \,
\sum_{k=1}^n \gamma_k=0, \, \, o_k \, \gamma_k=0, \, \, 1 \le k \le m
\rangle.
$

\begin{example} 
\label{assumptioniii)example0}
\normalfont
Let $\Sigma$ be the smooth orbisurface with underlying space $(\R/\Z)^2$ and with
exactly one cone point, say $p_1$, of order $o_1=2$. Let $\widetilde{\gamma}$ be a boundary loop
of  a small disk containing  $p_1$, and let $\gamma=[\widetilde{\gamma}]$. 
Let $\alpha, \, \beta$ define the standard basis of loops
corresponding to $(\R/\Z)^2$, that is, $\alpha,\, \beta$ are a
basis of the quotient free first orbifold homology group of $\Sigma$. For $x_0 \in \Sigma$,
$
\pi_1^{\textup{orb}}(\Sigma,\,x_0)=\langle \alpha, \, \beta, \, \gamma \, | \,
[\alpha, \, \beta]=\gamma, \, \gamma^2=1 \rangle,
$
and
$
\op{H}_1^{\textup{orb}}(\Sigma,\, \Z)=\langle \alpha, \, \beta, \, \gamma \, | \,
\gamma=1, \, 2 \, \gamma=0 \rangle \simeq \langle \alpha, \, \beta \rangle.
$
 \end{example}

Fix an $m$\--dimensional vector of strictly positive integers $\vec{o}$ and
a non\--negative integer $g$. Let
$\mathcal{M}\textup{S}^{\vec{o}}_m$ be the set of matrices
$B \in \textup{GL}(m, \, \Z)$ such that $B \cdot \vec{o}=\vec{o}\}$.
Define the quotient
$
T^{2g+m}_{(g;\, \vec{o})}/\mathcal{G}_{(g, \, \vec{o})}
$
 where
$
T^{2g+n}_{(g;\, \vec{o})}
$
is the set of points $(t_i)_{i=1}^{2g+m} \in T^{2g+m}$ such that 
the order of $t_i$ is  a multiple of  $o_i$ for all  $2g +1\le i \le m$,
and
$$
\mathcal{G}_{(g; \, \vec{o})}:=\Big\{  
\left( \begin{array}{cc}
A & \textup{0}  \\
C & D 
\end{array} \right) \in \textup{GL}(2g+m, \, \Z)
 \, | \, A \in \textup{Sp}(2g, \, \Z), \, \,  D \in \mathcal{M}\textup{S}^{\vec{o}}_m \Big\}.
$$
We call $
T^{2g+m}_{(g;\, \vec{o})}/\mathcal{G}_{(g, \, \vec{o})}
$ the \emph{Fuchsian signature} space associated to $(g;\, \vec{o})$
 \\
 \\
 In the case of our compact, connected symplectic manifold $(M,\omega)$ with an effective symplectic 
 $T$\--action with $\dim T$\--dimensional symplectic
 $T$\--orbits we have the following ingredients: let $(g; \, \vec{o}) \in \Z^{1+m}$ be the Fuchsian signature
of  $M/T$; let $\{\gamma_k\}_{k=1}^m$ be a basis of small loops around the cone points
$p_1,\ldots,p_n$ of $M/T$, viewed as an orbifold; let $\{\alpha_1, \, \beta_1,\ldots, \alpha_g,\beta_g\}$ be a symplectic basis of a \emph{free} subgroup 
$F$ of the orbifold homology group $\op{H}_1^{\textup{orb}}(M/T, \, \Z)$,
whose direct sum with the \emph{torsion} subgroup is $\op{H}_1^{\textup{orb}}(M/T, \, \Z)$;
let $\mu_{\textup{h}}$ be the homomomorphism induced on homology by
the monodromy homomorphism $\mu$ of the connection $\Omega$ of
symplectic orthogonal complements to the tangent spaces to the $T$\--orbits. 
The \emph{monodromy invariant of $(M, \omega, \, T)$} is the  
$\mathcal{G}_{(g; \, \vec{o})}$\--orbit
\begin{eqnarray} \label{mono}
 \mathcal{G}_{(g, \, \vec{o})} \cdot ((\mu_{\textup{h}}(\alpha_i), \, \mu_{\textup{h}}(\beta_i))_{i=1}^{g}, \, (\mu_{\textup{h}}(\gamma_k))_{k=1}^m)  \in T^{2g+m}_{(g;\, \vec{o})}/\mathcal{G}_{(g, \, \vec{o})}.
 \end{eqnarray}
Even though (\ref{mono}) depends on choices, one
can show that it is  well\--defined.

\subsubsection{Invariants of $M=S^2 \times_{\Z/2\,\Z} (\R/\Z)^2$} The invariants are: the non\--degenerate antisymmetric bilinear form
$
\omega^{\mathbb{R}^2}=\left( \begin{array}{cc}
0 & 1  \\
-1 & 0 
\end{array} \right)
$;  2) The Fuchsian signature $(g;\, \vec{o})=(0;\,2,\,2)$ of the orbit space $M/\T^2$; 3)
The symplectic area of  $S^2/(\mathbb{Z}/2\,\mathbb{Z})$: $1$ (half of the area of $S^2$); 4)
The monodromy invariant: 
$$
 \mathcal{G}_{(0; \, 2, \,2)} \cdot (\mu_{\textup{h}}(\gamma_1), \, \mu_{\textup{h}}(\gamma_2))=
  \Big\{ \left( \begin{array}{cc}
1 & 0  \\
0 & 1 
\end{array} \right), \, 
  \left( \begin{array}{cc}
0 & 1  \\
1 & 0 
\end{array} \right)
                 \Big \} \cdot ([1/2, \, 0],\, [1/2, \, 0]).$$
Here the $\gamma_1,\gamma_2$ are small loops around the poles of $S^2$. Then $M/T=S^2 / (\Z/2\,\Z)$,
$\pi^{\textup{orb}}_1(M/T,\,p_0)=\langle \gamma_1\, |\, \gamma_1^2=1 \rangle \simeq \Z/2\,\Z$, and
$\mu \colon \langle \gamma_1\, |\, \gamma_1^2=1 \rangle  \to T= (\R/\Z)^2$ is 
$
\mu(\gamma_1)=[1/2, \,0].
$
We have a $T$\--equivariant symplectomorphism
$$
\widetilde{M/T} \times_{\pi^{\textup{orb}}_1(M/T,\,p_0)} T 
=\widetilde{S^2 / (\Z/2\,\Z)} \times_{\pi^{\textup{orb}}_1(S^2 / (\Z/2\,\Z),\,p_0)} (\R/\Z)^2 
\simeq M.
$$

\section{Symplectic $2$\--torus actions on $4$\--manifolds} \label{four:sec}

\normalfont

Our only assumption now is that $T$ is  a $2$\--dimensional torus and $(M,\omega)$ is  a 
compact, connected $4$\--dimensional symplectic manifold on which $T$ acts effectively and 
symplectically.
The following is a simplified version of the main result of \cite{Pe}; readers can consult \cite[Theorem~8.2.1]{Pe} for the most informative and
complete version of the statement.

\begin{theorem}[\cite{Pe}] \label{main:thm}
If $(M,\omega)$ is compact, connected, symplectic $4$\--manifold equipped with an effective symplectic action of
a $2$\--torus  $T$ then one and only one of the following four cases occurs.
\begin{itemize}
\item[{1)}]
$(M, \, \omega)$ is a symplectic toric $4$\--manifold. 
\item[{2)}]
$(M, \, \omega)$ is equivariantly symplectomorphic to $(\R/\Z)^2 \times S^2$, where on  $(\R/\Z)^2 \times S^2$ 
the symplectic form is a product form.
The action of the $2$\--torus is: one circle acts on the first circle of $(\R/\Z)^2$ by translations,
while the other circle acts on $S^2$ by rotations about the vertical axes.
\item[{3)}]
$(M, \, \omega)$ is equivariantly symplectomorphic to $(T \times \mathfrak{t}^*)/\iota(P)$,
where
$T \times \got{t}^*$ is equipped with the standard cotangent bundle form and the standard $T$\--action on left factor of $T \times \mathfrak{t}^*$, both of which descend to $(T\times \got{t}^*)/\iota(P)$,
$P \subset \mathfrak{t}^*$ is a discrete cocompact subgroup,  and
$\iota(P) \le T \times \mathfrak{t}^*$ is a discrete, cocompact, subgroup for
a non\--standard group structure on $(T\times \got{t}^*)/\iota(P)$ defined in terms of
the Chern class (a certain bilinear form $c \colon\mathfrak{t}^* \times \mathfrak{t}^* \to \mathfrak{t}$), and
the holonomy invariant
$
[\tau \colon P \to T]_{\textup{exp}(\mathcal{A})} \in \op{Hom}_c(P,\, T)/{\rm exp}(\mathcal{A}).
$
\item[{4)}]
$(M, \, \omega)$ is equivariantly symplectomorphic to
$\widetilde{\Sigma} \times_{\pi^{\textup{orb}}_1(\Sigma, \, p_0)}T
$
where the symplectic form and $T$\--action are induced by the product ones,
$\Sigma$ is a good orbisurface,
and the orbifold fundamental group $\pi^{\textup{orb}}_1(\Sigma, \, p_0)$ acts on $\widetilde{\Sigma} \times T$ diagonally,
where the action of $\pi_1^{\scriptop{orb}}(\Sigma,\,p_0)$ on $T$ is by means
of any homomorphism $\mu \colon \pi_1^{\scriptop{orb}}(\Sigma) \to T$.
\end{itemize}
\end{theorem}

The proof of Theorem~\ref{main:thm} uses as stepping stones the symplectic case treated in Section~\ref{symp:sec},
and the coisotropic case treated in Section~\ref{cos:sec}.
Case 4) comes from Theorem~\ref{t3} and  Theorem~\ref{t4}; that is, the symplectic invariants and the
construction of this case comes from these results. The  fundamental observation to use the results in the aforementioned sections in order to prove Theorem~\ref{main:thm} is that under the 
assumptions of Theorem~\ref{main:thm}, 
there are precisely two possibilities: a) either all of the $T$\--orbits are symplectic $2$\--tori; or
b)  all of the $2$\--dimensional $T$\--orbits are Lagrangian $2$\--tori. A significant part of 
the proof of Theorem~\ref{main:thm} 
consists of unfolding item b) above into items 1), 2), 3) in the statement
of Theorem~\ref{main:thm}. Notice that item 1)  
is classified in terms of the Delzant polytope, which is the only symplectic
invariant in that case, see Remark~\ref{dez}, part (c).

\section{Comments on the Proofs} \label{proof:sec}

In order to give an idea of the type of methods involved in the study of non Hamiltonian
symplectic actions, we next give a outline of the proof of one of the results we have
stated earlier: Theorem~\ref{t1} (and readers can consult the references given througout
the paper for the proofs of the other results, and for more detailed and general versions
of them). For the notation and ingredients which we use next
consult Section~\ref{t1sec}. We have divided the outline of the proof (in \cite{DuPe}) 
into three steps:

{Step 1.}  \emph{The orbit space}. First one shows that the orbit space $M/T$ is a polyhedral $\mathfrak{l}^*$\--parallel space, where $\mathfrak{l}:=\textup{ker}(\omega^{\mathfrak{t}})$. Then one shows that, as 
$\mathfrak{l}^*$\--parallel spaces, there is an isomorphism 
$
M/T \simeq \Delta \times S. 
$
Here $\Delta$ is a Delzant polytope, and $S$ is a  torus. 
The polytope $\Delta$ encodes the Hamiltonian part of the $T$\--action, while
the torus $S$ encodes the free part of the action, see~Section~\ref{362:sec}.

{Step 2.}  \emph{A nice connection}. We prove that there exists an admissible connection for the principal
torus bundle $\pi \colon M_{\textup{reg}} \to M_{\textup{reg}}/T$ (i.e. a ``nice" linear assignment of
lifts) 
$
\xi \in \mathfrak{l}^* \mapsto L_{\xi} \in \mathcal{X}^{\infty}(M_{\textup{reg}}),
$
that has ``simple" Lie brackets $[L_{\xi}, \,L_{\eta}]$ in the sense that 
$[L_{\xi}, \,L_{\eta}]=c(\xi,\,\eta)_M$ if $\xi,\,\eta \in N:=(\mathfrak{l}/\mathfrak{t}_{\textup{h}})^*$ 
($c$ encodes Chern class of $\pi \colon M_{\textup{reg}} \to M_{\textup{reg}}/T$), and 
$[L_{\xi}, \,L_{\eta}]=0$ otherwise,
as well as ``simple" symplectic pairings $\omega(L_{\xi}, \, L_{\eta})$. 
The vector fields $L_{\xi}$, $\xi \in N$, have smooth extensions to $M$ but
 are singular on $M \setminus M_{\textup{reg}}$ (see~Section~\ref{362:sec}).

{Step 3.}  \emph{Distribution by symplectic toric manifolds}. 
We define the distribution 
$
D_x:=\textup{span}\{\, L_{\eta}(x),\,\,Y_M(x) \,\,\,| \,\,\, Y \in \mathfrak{t}_{\textup{h}}, \,\, \eta \in C\,\},\,\,\, 
x \in M_{\scriptop{reg}},
$
on $M_{\scriptop{reg}}$ and prove that it is integrable, where $C \oplus N = \mathfrak{l}^*$. 
The integral manifolds to this distribution are $T$\--equivariantly symplectomorphic to
$(M_{\textup{h}},\,\omega_{\textup{h}}, \,T_{\textup{h}})$.
 Next we give a natural construction of the groups $H$, $G$, and $G \times_H M_{\textup{h}}$ from
the connection $\xi \mapsto L_{\xi}$. The group
$G$ involves the Chern class of $M_{\scriptop{reg}} \to M_{\scriptop{reg}}/T$
and $H$ involves the holonomy of $\xi \mapsto L_{\xi}$. We gave the precise formulas in Section~\ref{Gdef}
and Section~\ref{hdef:sec}.
 Then we show that the following map is a $T$\--equivariant
symplectomorphism from $G \times_H M_{\textup{h}}$ to $M$:
$
((t,\,\xi),\,x) \mapsto t \cdot \textup{e}^{L_{\xi}}(x). 
$
This concludes the sketch of proof of Theorem~\ref{t1}.

\medskip

This paper is dedicated to the memory of Professor Johannes J. Duistermaat (1942--2010). 
The article \cite{GuPeVNWe} outlines some of his  contributions (see also \cite[Section 2.4]{pevn11}). 
The author is supported by  NSF CAREER DMS-1518420, and   
thanks Joseph Palmer for comments on a preliminary version.

\smallskip\noindent
\noindent
{{\'A}lvaro Pelayo}\\
University of California, San Diego\\ 
Department of Mathematics \\
9500 Gilman Dr \#0112\\
La Jolla, CA 92093-0112, USA \\
{\em E\--mail}: \texttt{alpelayo@math.ucsd.edu} \\

\end{document}